\documentclass[a4paper]{article}

\usepackage{amsmath,amsthm,amssymb,mathtools}
\allowdisplaybreaks

\usepackage[margin=2cm]{geometry}
\usepackage{url}

\usepackage{xcolor}

\usepackage{hyperref}

\numberwithin{equation}{section}
\newtheorem{theorem}{Theorem}[section]

\newtheorem{lemma}{Lemma}[section]

\newtheorem{remark}{Remark}[section]

\def\RR{\mathbb{R}}
\def\NN{\mathbb{N}}
\def\Cspace#1#2{C^{#1}(#2)}

\def\Lspace#1#2{L^{#1}(#2)}
\def\Wspace#1#2#3{W^{#1,#2}(#3)}

\def\tmax{T_{\text{max}}}
\def\cGN{C_{\text{GN}}}

\def\intO{\int_\Omega}
\def\ngrad#1{\lvert\nabla #1 \rvert}
\def\dert{\frac{\text{d}}{\text{d}t}}
\def\ds{\text{d}s}
\def\measDom{\lvert\Omega\rvert}
\def\cY#1#2{C_{\text{Y}}\left(#1,#2\right)}

\def\partialn#1{\frac{\partial #1}{\partial\nu}}

\usepackage{constants}

\newconstantfamily{const}{
	symbol=c, 
}
\newconstantfamily{tau}{
	symbol=\tau, reset={section}
}
\newconstantfamily{eps}{
	symbol=\varepsilon, reset={section}
}

\usepackage{hyperref}

\title{Tumor-immune cell interactions by a fully parabolic chemotaxis model with logistic source}
\author{Rafael D\'iaz Fuentes}
\date{}

\begin{document}
	
	\maketitle
	\centerline{$^{\natural}$Department of Mathematics and Computer Science, University of Cagliari, Cagliari, Italy} 
	\centerline{({\tt rafael.diazfuentes@unica.it}) }
	
	\begin{abstract} 
		This work studies the existence of classical solutions to a class of chemotaxis
		systems reading
		\begin{equation*}
			\begin{dcases}
				u_t = \Delta u-\chi \nabla\cdot(u \nabla v) + \mu_1 u^k -\mu_2 u^{k+1},\hspace*{0.5cm} & \text{in} \; \Omega\times(0,\tmax),\\
				v_t= \Delta v+\alpha w-\beta v-\gamma u v, & \text{in} \; \Omega\times(0,\tmax),\\
				w_t= \Delta w-\delta u w+ \mu_3 w(1-w), & \text{in} \; \Omega\times(0,\tmax),\\
				\partialn{u}=\partialn{v}=\partialn{w}=0, & \text{on} \; \partial\Omega\times(0,\tmax),\\
				u(x,0)=u_0(x), \quad v(x,0)= v_0(x), \quad w(x,0)= w_0(x), & x\in\overline{\Omega},
			\end{dcases}
		\end{equation*}
		that model interactions between tumor (i.e., $w$) and immune cells (i.e., $u$) with a logistic-type source term $\mu_1 u^k - \mu_2 u^{k+1}$, $k\geq1$, also in presence of a chemical signal (i.e., $v$). The model parameters $\chi, \mu_1,\mu_2, \mu_3, \alpha, \beta, \gamma$, and $\delta$ are all positive. The value $\tmax$ indicates the maximum instant of time up to which solutions are defined. Our focus is on examining the global existence in a bounded domain $\Omega\subset \RR^n, n \geq 3$, under Neumann boundary conditions. We distinguish between two scenarios: $k>1$ and $k=1$. The first case allows to prove boundedness under smaller assumptions relying only on the model parameters instead of on the initial data, while the second case requires an extra condition relating the parameters $\chi, \mu_2$, $n$, and the initial data $\lVert v_0 \rVert_{L^\infty(\Omega)}$. This model can be seen as an extension of those previously examined in \cite{LankeitWang17} and \cite{GnanasekaranCDFN24}, being the former a system with only two equations and the latter the same model without logistic. 
		
		\medskip
		\noindent\textbf{Keywords:} chemotaxis systems, classical solutions, global existence, tumor-immune cell interaction.
		
		\noindent\textbf{2020 Mathematics Subject Classification:} 35A01, 35A09, 35B30, 35B45, 92C17
	\end{abstract}

	\section{Introduction}
	
	This work studies the dynamics between tumor cells and immune cells influenced by the presence of a chemi\-cal signal released by the latter. The interplay between these units has been broadly studied as a chemotaxis model \cite{BellomoBTW15,Painter19}. Chemotaxis, as a model where units (microorganisms as bacteria or organisms as animals \cite{TelloWrzosek16}) move in response to chemical gradients, was first modeled by Keller and Segel considering mainly two dynamics. A first one is the so-called \emph{production model} \cite{KellerSegel70prod}, where a population of cells $u=u(x,t)$ (e.g., amoebae in that first study), $x \in \Omega$ a bounded region of $\RR^2$, $t\in\RR$, releases a chemical substance $v=v(x,t)$ (in that case the acrasin), related to the simplified mathematical model:
	\begin{equation} \label{eq:KSprod}
		\begin{dcases}
			u_t = \Delta u - \chi \nabla \cdot (u \nabla v), & \text{in} \ \Omega\times(0,\tmax) \\
			v_t = \Delta v - v + u, & \text{in} \ \Omega\times(0,\tmax) \\ 
			\partialn{u} = \partialn{v} = 0, & \text{on} \; \partial\Omega\times(0,\tmax),\\
			u(x,0)=u_0(x), \quad v(x,0)= v_0(x),  & x\in\overline{\Omega},
		\end{dcases}
	\end{equation}  
	with $\chi$ the {\em chemotactic coefficient} measuring the strength of chemotaxis, i.e., the term $-\chi\nabla\cdot(u\nabla v)$ represents the movement of the cells towards the chemical signal. The vector $\nu$ is the outward unit normal on $\partial \Omega$. The value $\tmax \in (0,\infty]$ indicates the maximum instant of time up to which solutions are defined. All the models discussed within this work, consider homogeneous Neumann boundary conditions and therefore, from now on, that condition (e.g., the third line in \eqref{eq:KSprod}) will be omitted.
	
	The second is the \emph{consumption model} \cite{KellerSegel71cons}, reading in a simplified form:
	\begin{equation} \label{eq:KScons}
		\begin{dcases}
			u_t = \Delta u - \chi \nabla \cdot (u \nabla v), & \text{in} \ \Omega\times(0,\tmax) \\
			v_t = \Delta v - uv, & \text{in} \ \Omega\times(0,\tmax) \\ 
			u(x,0)=u_0(x), \quad v(x,0)= v_0(x),  & x\in\overline{\Omega},
		\end{dcases}
	\end{equation}
	where the population of bacteria $u(x,t)$ moves toward higher concentration of oxygen $v(x,t)$ depleting it. For more details on these models and the main theoretical techniques to tackle their studies the reader may address the surveys and papers \cite{BellomoBTW15,HorstmannSurvey03, HorstmannSurvey04, Horstmann2005, Tao11, TaoWinkler12A, TaoWinkler12C, BaghaeiKelghati17, Painter19, LankeitWinklerSurvey20, LankeitWinklerSurvey23}. Those works deal with analyses concerning local and/or global solutions for models defined on bounded and smooth domains $\Omega$ of $\RR^n$, for $n\geq 1$; specifically some aspects regarding the possibility of blow-up in finite time. 
	
	Just to summarize some of them, regarding the model \eqref{eq:KSprod}, as mentioned in \cite{TaoWinkler12A} and the references cited there, for $n=1$ all the solutions are global in time and bounded; meanwhile, for $n=2$ then it is required that 
	\begin{math}
		\lVert u_0 \rVert_{\Lspace{1}{\Omega}} < 4\pi.
	\end{math}
	Otherwise, for $n\geq 3$, global and bounded solutions exist provided that 
	\begin{math}
		\lVert u_0 \rVert_{\Lspace{\frac{n}{2}+\delta}{\Omega}} < \varepsilon
	\end{math}
	and 
	\begin{math}
		\lVert \nabla{v_0} \rVert_{\Lspace{n+\delta}{\Omega}} < \varepsilon,
	\end{math}
	for some $\varepsilon>0$ and any $\delta>0$ \cite{Winkler10}. However, blow up either in finite or infinite time may occur for $n=2$ if 
	\begin{math}
		\lVert u_0 \rVert_{\Lspace{1}{\Omega}} > 4\pi.
	\end{math}
	Also for $n\geq3$, the solution may be unbounded either in finite or infinite time when $\Omega$ is a ball. 
	
	On the other hand, for the consumption model, in \cite{Tao11} it was proved that the system \eqref{eq:KScons} possesses a unique classical global uniformly-in-time bounded solution if
	\begin{math}
		0 < \chi \leq \frac{1}{6(n+1) \lVert v_0 \rVert_{L^\infty(\Omega)}},
	\end{math}
	for a smooth bounded domain $\Omega \subset \RR^n$ and $n\geq 2$. We should notice the required relation between the parameter $\chi$, the initial condition $v_0(x)$, and the dimension $n$. That restriction was sharpened by \cite{BaghaeiKelghati17}, always linking those three values. In case of arbitrary large initial data, global weak solutions are obtained in smooth bounded convex domains $\Omega \subset \RR^3$ \cite{TaoWinkler12C}.
	
	The aggregation processes mentioned above are essentially connected to an uncontrolled increasing of the cells distribution; this is typical for production models \eqref{eq:KSprod}, where the intensity of attraction actions are higher. In this sense, the presence of an external source behaving as $u - u^k$, $k>1$ (i.e., logistic-type), in the equation for $u$ provides dampening effects possibly limiting the presence of blow up; despite that, for $k$ close to 1, even this presence is not sufficient to prevent gathering formation \cite{TelloWinkler07,Winkler18}. 
	
	Conversely, for consumption models \eqref{eq:KScons} (which are of our interest in this paper), the situation is far to be similar. Indeed, even though this is a challenging problem, it is conjectured that such mechanisms do not present explosion phenomena. As indicated above, globality and boundedness are in any case achieved provided smallness assumptions on $\chi \lVert v_0 \rVert_{\Lspace{\infty}{\Omega}}$. As a consequence, it appears natural introducing logistic sources not in order to prevent blow up but to ensure boundedness even for large values of $\chi \lVert v_0 \rVert_{\Lspace{\infty}{\Omega}}$. 
	This choice can be traced in the model studied by Lankeit and Wang \cite{LankeitWang17} 
	\begin{equation} \label{eq:modelLankeit}
		\begin{dcases}
			u_t = \Delta u - \chi \nabla\cdot(u \nabla v) + \mu_1 u -\mu_2 u^{2},\hspace*{0.5cm} & \text{in} \; \Omega\times(0,\tmax),\\
			v_t= \Delta v- u v, & \text{in} \; \Omega\times(0,\tmax),\\
			u(x,0)=u_0(x), \quad v(x,0)= v_0(x), & x\in\overline{\Omega},
		\end{dcases}
	\end{equation}
	with $\Omega \subset \RR^n$ a smooth bounded domain, $n\in\NN$.
	The presence of such an external source of logistic-type benefits the boundedness of the global solution provided a restriction on the problem parameter $\mu_2$ reading
	\begin{equation}\label{eq:muLankeit}
		\mu_2 > K_1(n) \lVert \chi v_0 \rVert_{L^\infty(\Omega)}^{\frac1n} + K_2(n) \lVert \chi v_0 \rVert_{L^\infty(\Omega)}^{2n},
	\end{equation}
	for some positive constants $K_1 = K_1(n)$ and $K_2 = K_2(n)$.
	
	Several extensions to those models have been proposed including as well a third equation in model relating a prey and a predator \cite{TelloWrzosek16}. In this context, the predator moves toward the gradient of concentration of some chemical released by prey instead of moving directly toward the higher density of prey \cite{WangWang20}. The dynamics of prey and predator consider both, the production and consumption models; in the sense that an unit produce another, while it is consumed by a third one. 
	Let us consider as a reference point the model
	\begin{equation} \label{eq:modelGen}
		\begin{dcases}
			u_t = \Delta u - \chi \nabla \cdot (u \nabla v) + g(u) & \text{in} \; \Omega\times(0,\tmax), \\
			v_t = \Delta v + h_1(u,v,w) & \text{in} \; \Omega\times(0,\tmax), \\
			w_t = \Delta w + h_2(u,w) & \text{in} \; \Omega\times(0,\tmax),
		\end{dcases}
	\end{equation}
	and the respective initial conditions for each function. 
	
	In \cite{HuTao20}, Hu and Tao identify a critical mass of lymphocytes necessary for the coexistence of tumor and immune cells, highlighting the importance of immune response in tumor dynamics. The proposed model describes the interaction of two cells with one chemical secreted by the tumor cells, where $u=u(x,t)$ represents the density of lymphocytes, $v=v(x,t)$ describes the chemical signal secreted by tumor cells and $w=w(x,t)$ denotes the density of tumor cells. In this case, $x\in\Omega\subset\RR^2$, a bounded convex domain with smooth boundary. The term $-\chi\nabla\cdot(u\nabla v)$ in the first equation in \eqref{eq:modelGen} represents the movement of the lymphocytes towards the chemical signal, and the external source is ignored, i.e., $g(u)=0$. The sources for the second and third equations are chosen as $h_1(u,v,w)= w - v - uv$, which combines the production of chemical signal and the consumption by the immune cells, and $h_2(u,w) = -uw + w(1-w)$ combining the consumption of the tumor cells by the immune cells with a natural logistic growth term of the tumor cells. The system is proven to admit a unique global uniformly-in-time bounded classical solution for all nonnegative initial conditions $(u_0, v_0, w_0) \in {{C^{0}}(\overline{\Omega})} \times {{W^{1,\infty}}(\Omega)} \times {{ W^{1,\infty}}(\Omega)}$, with $u_0 \not \equiv 0$, $v_0>0$, and $w_0\not \equiv 0$.
	
	The lack of parameters scaling each term in the equations for their choices of $h_1(u,v,w)$ and $h_2(u,v)$, loses the sense of dimensionality of the model, as observed in \cite{TelloWrzosek16}. That was a motivation in \cite{GnanasekaranCDFN24} to propose a model with  $h_1(u,v,w)= \alpha w - \beta v - \gamma uv$ and $h_2(u,w) = -\delta uw + \mu w(1-w)$. The parameters $\chi$, $\alpha$, $\beta$, $\gamma$, $\delta$, $\mu$ are assumed positive constants. The constant $\alpha$ is the production rate of the chemical by the tumor cells, $\beta$ is the decay rate of the chemical and $\gamma$ is the consumption rate of the chemical by the lymphocytes. The parameter $\delta$ models the destruction rate of the tumor cells by the lymphocytes. Lastly, $\mu$ is the growth coefficient for the cancer cells. 
	
	In that work, it is considered $\Omega\subset\RR^n$ a bounded domain, with $n\geq3$. The studied fully parabolic system admit a unique nonnegative global classical solution uniformly bounded-in-time, provided 
	\begin{equation} \label{eq:condGnana}
		\max\left\{\frac{\alpha}{\beta}, \frac{\alpha}{\beta}\|w_0\|_{\Lspace{\infty}{\Omega}}, \|v_0\|_{\Lspace{\infty}{\Omega}} \right\}<\frac{\pi}{\chi}\sqrt{\frac{2}{n}}.
	\end{equation}
	
	A natural following up is to consider an external source for the first equation $g(u) = \mu_1 u^k - \mu_2 u^{k+1}$, for $k\geq 1$ and positive $\mu_1$ and $\mu_2$. Let us consider the model, extending the proposals in \cite{LankeitWang17,GnanasekaranCDFN24}, describing the tumor-immune cell interactions by
	\begin{equation}\label{eq:model}
		\begin{dcases}
			u_t = \Delta u-\chi \nabla\cdot(u \nabla v) + \mu_1 u^k -\mu_2 u^{k+1},\hspace*{0.5cm} & \text{in} \; \Omega\times(0,\tmax),\\
			v_t= \Delta v+\alpha w-\beta v-\gamma u v, & \text{in} \; \Omega\times(0,\tmax),\\
			w_t= \Delta w-\delta u w+ \mu_3 w(1-w), & \text{in} \; \Omega\times(0,\tmax),\\
			u(x,0)=u_0(x), \quad v(x,0)= v_0(x), \quad w(x,0)= w_0(x), & x\in\overline{\Omega},
		\end{dcases}
	\end{equation}
	with $\chi,\alpha,\beta,\gamma,\delta$, and $\mu_3$ positive constants. The considered domain $\Omega\subset\RR^n$, with $n\geq 3$, is bounded and with smooth boundary $\partial \Omega$. 
	
	The logistic source term for lymphocytes $g(u,v) = \mu_1 u^k - \mu_2 u^{k+1}$ captures the dynamics of lymphocyte population growth and proliferation reflecting the self-regulating nature of the immune response in the presence of a tumor \cite{Wilkie17}, \cite{Mahlbacher19}. As the number of lymphocytes increases, the rate of growth slows down, eventually leading to a plateau in the population size. To the best of our knowledge, there are no studies carried out for the system \eqref{eq:model} with logistic source of that form.
	
	\subsection{Main result and structure of the paper}

	In order to formally present the aforementioned analysis, let us suppose that initial values $u_0$, $v_0$ and $w_0$ satisfy
	\begin{equation}\label{eq:initCond}
		\begin{dcases}
			u_0\in\Cspace{0}{\overline{\Omega}},\quad \mbox{with} \quad u_0 \geq 0\quad \mbox{in}\: \overline{\Omega},\\
			v_0, w_0\in \Wspace{1}{q}{\Omega},\quad \mbox{for some}\,\, q>n,\quad \mbox{with} \quad v_0, w_0 \geq 0\quad\mbox{in}\: \Omega.
		\end{dcases}
	\end{equation}
	
	With these initial conditions we are able to state what follows.
	
	\begin{theorem} \label{theo:globalexist}
		Let $\Omega\subset\RR^n$, with $n\geq3$, be a bounded domain with smooth boundary and let $\chi,\alpha,\beta,\delta,\gamma, \mu_1, \mu_2$, and $\mu_3$ be positive constants. Let $u_0, v_0$, and $w_0$ be initial conditions to \eqref{eq:model} fulfilling the hypotheses in \eqref{eq:initCond}. For $M_2 :=\max\Big\{1, \|w_0\|_{\Lspace{\infty}{\Omega}}\Big\}$ and $M_3:=\max \Big \{\frac{\alpha}{\beta}M_2, \|v_0\|_{\Lspace{\infty}{\Omega}}\Big\}$, let us assume that those parameters satisfy for $k > 1$
		\begin{equation} \label{eq:theocond}
			\delta > \mu_3 \quad \text{and} \quad \alpha < 2 \min\left\{\beta+\gamma,\delta-\mu_3\right\}
		\end{equation} 
		with an extra requirement for $k=1$ reading
		\begin{equation} \label{eq:muRaf}
			\mu_2 > \ A_1(n) \chi^{2+\frac4n} M_3^\frac4n
			+ A_2(n) \gamma^{\frac{n+2}{2}}  M_3^{n} + A_3(n) \delta^{\frac{n+2}{2}} M_2^{n},
		\end{equation} 
		for certain computable constants $A_1(n), A_2(n)$, and $A_3(n)$.
		
		Then, problem \eqref{eq:model} admits a unique classical solution, global and uniformly bounded in time.
	\end{theorem}

	\begin{remark}
		It is worth realizing that for $k>1$, large values of 
		\begin{math}
			\chi \lVert v_0 \rVert_{L^\infty(\Omega)}
		\end{math}
		are allowed and the restrictions are only on the problem parameters $\alpha, \beta, \delta, \gamma$, and $\mu_3$.
	\end{remark}
	
	\begin{remark} \label{rem:comparison}
		If we compare the condition in \eqref{eq:muLankeit} with \eqref{eq:muRaf}, considering model \eqref{eq:model} as a generalization for $k=1$ of model \eqref{eq:modelLankeit}, it becomes evident that adding that equation for $w$ does not change the type of constraint on the data; it only adds an additional term associated with $w$. However, for $k>1$ the benefit of the constraint that does not depend on 
		\begin{math}
			\lVert v_0 \rVert_{L^\infty(\Omega)}
		\end{math}
		and
		\begin{math}
			\lVert w_0 \rVert_{L^\infty(\Omega)}
		\end{math}
		is noticeable.
		
		On the other hand, comparing \eqref{eq:condGnana} with \eqref{eq:muRaf}, we can see how the extension of the model proposed in \cite{GnanasekaranCDFN24} to model \eqref{eq:model} has the advantage of eliminating the constraint for the value
		\begin{math}
			\lVert w_0 \rVert_{L^\infty(\Omega)}.
		\end{math}
	\end{remark}
	
	This theorem is proved as a consequence of some lemmas and quoted inequalities. To ease the argumentation of its validity the paper is organized as follows. 
	Section \ref{sec:Prelim} presents some known results from the literature which are used during the proofs in this paper. The next Section \ref{sec:localExist} proves the local existence of solutions $(u,v,w)$ for the considered model. These solutions are positive and their boundedness in some Sobolev spaces is proved. Moreover, an extensibility criterion is given which supports the globality of the existing local solution. Finally, Section \ref{sec:aprioriest} proves the main result already stated in this section (i.e., Theorem \ref{theo:globalexist}) as a consequence of some {\em a priori} inequalities and the results in Section \ref{sec:localExist}. To conclude this section, some thoughts are devoted to further refining the arguments in {\em Remark \ref{rem:comparison}}. 
	
	
	Along the proofs many constants, denoted by $c_i, i\in\NN$, are explicitly computed, supporting the statement in Theorem \ref{theo:globalexist} to provide computable constants. However, their expressions may be cumbersome for the reader in some cases. At the same time, showing them makes it easier for the reader to follow the passages along the inequalities, hence the decision to include an Appendix \ref{sec:contants} where the explicit values of some constants are shown.

	\section{Some preliminary results}
	\label{sec:Prelim}
	
	We dedicate this section to present some lemmas supporting the main results within this work.
	
	With the Gagliardo-Nirenberg inequality, we refer to the version proposed in \cite{LiLankeit16}. 
	
	Let us present some inequalities from \cite[Lemma 2.1, Lemma 2.2]{LankeitWang17} in the following lemma. The term $D^2 \psi$ denotes the Hessian of $\psi$ and $\lvert D^2 \psi \rvert^2 = \sum_{i,j=1}^n \left(\frac{\partial^2 \psi}{\partial x_i x_j}\right)^2$.
	
	\begin{lemma} \label{lem:lankeit}
		
		\begin{enumerate}
			\item For any function $\psi\in C^2(\Omega)$:
			\begin{equation} \label{ineq:LapHessian}
				\lvert\Delta\psi\rvert^2 \leq n \lvert D^2 \psi \rvert^2
			\end{equation}
			throughout $\Omega$.
			\item Let $p\in[1,\infty)$. For every $\eta>0$ there is $C_\eta>0$ such that every function $\psi\in C^2(\Bar{\Omega})$ with $\partialn \psi=0$ on $\partial \Omega$ satisfies
			\begin{equation}\label{ineq:boundaryGradIneq}
				\int_{\partial \Omega} \ngrad{\psi}^{2p-2}\ \partialn{\ngrad{\psi}^2} \leq \eta 	\intO \ngrad{\psi}^{2p-4}\ngrad{\ngrad{\psi}^2}^2+C_\eta \left(\intO \ngrad{\psi}^2\right)^p.
			\end{equation}
			\item Let $q \in [1, \infty)$. Then for any $\psi \in C^2(\overline{\Omega})$ satisfying $\psi \frac{\partial \psi}{\partial \nu} = 0$ on $\partial \Omega$, the inequality
			\begin{equation} \label{ineq:gradHessian}
				\lVert \nabla \psi\rVert_{L^{2(q+1)}(\Omega)}^{2(q+1)} \leq 2(4q^2 + n)\rVert \psi\rVert_{L^\infty(\Omega)}^2 \lVert\lvert\nabla \psi\rvert^{q-1} D^2 \psi\rVert_{L^2(\Omega)}^2
			\end{equation}
			holds.
		\end{enumerate}
	\end{lemma}
	
	It will turn out to be useful as well the following two known results which are presented here without proof.
	
	\begin{lemma}
		For any function  $\psi\in C^2(\Omega)$ it holds the identity
		\begin{equation} \label{eq:idLapNormGrad}
			\Delta \ngrad{\psi}^2 = 2\nabla \psi \cdot \nabla \Delta \psi + 2\lvert D^2 \psi \rvert^2.
		\end{equation}	
	\end{lemma}
	
	\begin{lemma}
		For any function $\psi\in C^2(\Omega)$ and $p>1$, it holds the inequality
		\begin{equation} \label{ineq:NgradNgrad}
			\ngrad{\ngrad{\psi}^p}^2 \leq p^2 \ngrad{\psi}^{2p-2} \lvert D^2 \psi \rvert^2.
		\end{equation}
	\end{lemma}

	The following Lemma presents some relations later used referring to the functions $u, v$, and $w$. Some of them interchange those functions, therefore the choice to use generic functions $\psi$ and $\varphi$ playing the role of solutions of the problem, i.e., satisfying the conditions obtained in {\em Lemma \ref{lem:localex}}.
	
	\begin{lemma} \label{lem:ineq}
		For any nonnegative functions $\psi,\varphi \in \Cspace{2}{\overline{\Omega}}$, any positive $\Cl[tau]{t:uYoung}, \Cl[tau]{t:gradvYoung}$, and any $p>1$, the following inequalities are valid:
		\begin{enumerate}
			\item for any $k>1$ there exists a computable positive constant $\Cl[const]{ce:uPgradv2}(\Cr{t:uYoung},\Cr{t:gradvYoung})$ such that 
			\begin{equation} \label{ineq:uPgradv2}
				\intO \psi^p \ngrad{\varphi}^2 \leq \Cr{t:uYoung} \intO \psi^{p+k} + \Cr{t:gradvYoung} \intO \ngrad{\varphi}^{2(p+1)} + \Cr{ce:uPgradv2}(\Cr{t:uYoung},\Cr{t:gradvYoung});
			\end{equation}  
			\item there exists a computable positive constant $\Cl[const]{ce:uPgradv2KO}(\Cr{t:uYoung})$ such that
			\begin{equation} \label{ineq:uPgradv2KO}
				\intO \psi^p \ngrad{\varphi}^2 \leq \Cr{t:uYoung} \intO \ngrad{\varphi}^{2(p+1)} + \Cr{ce:uPgradv2KO}(\Cr{t:uYoung}) \intO \psi^{p+1};
			\end{equation}
			\item for any $k>1$ there exists a computable positive constant $\Cl[const]{ce:u2gradvDPMD}(\Cr{t:uYoung},\Cr{t:gradvYoung})$ such that 
			\begin{equation} \label{ineq:u2gradvDPMD}
				\intO \psi^2 \ngrad{\varphi}^{2p-2} \leq \Cr{t:uYoung} \intO \psi^{p+k} + \Cr{t:gradvYoung} \intO \ngrad{\varphi}^{2(p+1)} + \Cr{ce:u2gradvDPMD}(\Cr{t:uYoung},\Cr{t:gradvYoung});
			\end{equation}
			\item there exists a computable positive constant $\Cl[const]{ce:u2gradvDPMDKO}(\Cr{t:uYoung})$ such that
			\begin{equation} \label{ineq:u2gradvDPMDKO}
				\intO \psi^2 \ngrad{\varphi}^{2p-2} \leq \Cr{t:uYoung} \intO \ngrad{\varphi}^{2(p+1)} + \Cr{ce:u2gradvDPMDKO}(\Cr{t:uYoung}) \intO \psi^{p+1};
			\end{equation}
			\item there exists a computable positive constant $\Cl[const]{ce:gradVgradW}(\Cr{t:uYoung},\Cr{t:gradvYoung})$ such that 
			\begin{equation} \label{ineq:gradVgradW}
				\intO \ngrad{\psi}^{2p-1} \ngrad{\varphi} \leq \Cr{t:uYoung} \intO \ngrad{\psi}^{2(p+1)} + \Cr{t:gradvYoung} \intO \ngrad{\varphi}^{2(p+1)} + \Cr{ce:gradVgradW}(\Cr{t:uYoung},\Cr{t:gradvYoung})
			\end{equation}
			\item there exists a computable positive constant $\Cl[const]{ce:gradW}(\Cr{t:uYoung})$ such that
			\begin{equation} \label{ineq:gradW}
				\intO \ngrad{\psi}^{2p} \leq \Cr{t:uYoung} \intO \ngrad{\ngrad{\psi}^p}^2 +  \Cr{ce:gradW}(\Cr{t:uYoung}).
			\end{equation}
		\end{enumerate}
		
		\begin{proof}
			Since the proof of some of the inequalities follows a similar line of reasoning, let us then show the proof for just a few of them. For instance, to prove the inequalities in \eqref{ineq:uPgradv2}, \eqref{ineq:u2gradvDPMD}, and \eqref{ineq:gradVgradW}, we can consider the first one of them. By the Young inequality
			\begin{equation*}
				\begin{split}
					\intO \psi^p \ngrad{\varphi}^2 \leq&\ \Cr{t:uYoung} \intO \psi^{p+k} + \tfrac{k}{p+k} \left(\tfrac{p}{\tau_1 (p+k)}\right)^{\frac{p}{k}} \intO \ngrad{\varphi}^{\frac{2(p+k)}{k}} \\
					\leq&\ \Cr{t:uYoung} \intO \psi^{p+k} + \Cr{t:gradvYoung} \intO \ngrad{\varphi}^{2(p+1)} + \Cr{ce:uPgradv2}(\Cr{t:uYoung},\Cr{t:gradvYoung}),
				\end{split}
			\end{equation*}
			with $\Cr{ce:uPgradv2}(\Cr{t:uYoung},\Cr{t:gradvYoung})$ as in \eqref{eqA:ceuPgradv2}.
			
			On the other hand, the inequalities \eqref{ineq:uPgradv2KO} and \eqref{ineq:u2gradvDPMDKO} are similarly derived. In particular, for \eqref{ineq:uPgradv2KO} we have
			\begin{equation*}
				\intO \psi^p \ngrad{\varphi}^2 \leq \Cr{t:uYoung} \intO \ngrad{\varphi}^{2(p+1)} + \Cr{ce:uPgradv2KO}(\Cr{t:uYoung}) \intO \psi^{p+1}
			\end{equation*}
			with $\Cr{ce:uPgradv2KO}(\Cr{t:uYoung})$ defined in \eqref{eqA:ceuPgradv2KO}.
			
			Lastly, regarding \eqref{ineq:gradW}, from the {\em AM-QM} inequality (i.e., $(a+b)^2 \leq 2(a^2+b^2$)) and the Gagliardo-Nirenberg inequality, we deduce
			\begin{equation} \label{eq:ineqgradWProof}
				\begin{split}
					\intO \ngrad{\psi}^{2p} &\leq\ 2 \cGN \lVert \nabla{\ngrad{\psi}^p}\rVert^{2\theta}_{\Lspace{2}{\Omega}} \lVert \ngrad{\psi}^p \rVert^{2(1-\theta)}_{\Lspace{\frac2p}{\Omega}} + 2\cGN \lVert \ngrad{\psi}^p \rVert^2_{\Lspace{\frac2p}{\Omega}} \\
					&\leq\ 2 \cGN \left(\intO \ngrad{\ngrad{\psi}^p}^2\right)^\theta \left(\intO \ngrad{\psi}^2\right)^{p(1-\theta)} + 2\cGN \left(\intO \ngrad{\psi}^2\right)^p, 
				\end{split}
			\end{equation}
			with $\theta=\frac{\frac{p}{2}-\frac12}{\frac{p}{2}-\frac12+\frac1n} \in (0,1)$. 
			
			For a reason which will be clear later, we use the bound $\intO \ngrad{\psi}^2 < C$ for certain constant $C>0$. Thus, the Young inequality applied to \eqref{eq:ineqgradWProof} leads to
			\begin{equation*}
				\begin{split}
					\intO \ngrad{\psi}^{2p} &\leq\ 2 \cGN C^{p(1-\theta)} \left(\intO \ngrad{\ngrad{\psi}^p}^2\right)^\theta + 2\cGN C^p \\
					&\leq\ \Cr{t:uYoung} \intO \ngrad{\ngrad{\psi}^p}^2 +  \Cr{ce:gradW}(\Cr{t:uYoung})
				\end{split}
			\end{equation*}
			with $\Cr{ce:gradW}(\Cr{t:uYoung})$ defined in \eqref{eqA:cegradW}. This complete the proof of \eqref{ineq:gradW}. 
		\end{proof}
	\end{lemma}

	\section{Local existence and extensibility criterion}
	\label{sec:localExist}
	
	At first, we state the local solvability of the system of equations \eqref{eq:model} and we prove some properties of the local solution $(u,v,w)$, as for instance the boundedness of $u$ in $L^1(\Omega)$ and of $v,w$ in $L^\infty(\Omega)$. An extensibility criterion is presented as support to derive the globality of the local solution obtained.
	
	\begin{lemma}[Local existence and positivity] \label{lem:localex}
		Suppose that $\Omega\subset \RR^n$, $n\geq 2,$ is a bounded domain with smooth boundary and $q>n$. Then for each nonnegative initial data satisfying \eqref{eq:initCond} and, for $k\geq 1$, $\chi,\alpha,\beta,\gamma,\delta,\mu_1, \mu_2$, and $\mu_3$ positive constants, there exists $\tmax\in (0,\infty]$ such that  the system \eqref{eq:model} admits a unique nonnegative solution $(u, v, w)$ belonging to
		\begin{align*}
			u&\in { C^{0}}\left(\overline{\Omega}\times\left.\left[0,\tmax\right.\right)\right)\cap { C^{2,1}}\left(\overline{\Omega}\times\left(0,\tmax\right)\right),  \\
			v, w&\in { C^{0}}\left(\overline{\Omega}\times\left.\left[0,\tmax\right.\right)\right)\cap { C^{2,1}}\left(\overline{\Omega}\times\left(0,\tmax\right)\right)\cap { L^{\infty}_{loc}}\left(\left.\left[0,\tmax\right.\right);\Wspace{1}{q}{\Omega}\right). 
		\end{align*}
		Furthermore, if $\tmax<\infty$, then 
		\begin{equation*} 
			\begin{split}
				\limsup_{t\to \tmax} \left(\lVert u(\cdot,t) \rVert_{\Lspace{\infty}{\Omega}} 
				+ \lVert v(\cdot,t) \rVert_{\Wspace{1}{q}{\Omega}} 
				+ \lVert w(\cdot,t) \rVert_{\Wspace{1}{q}{\Omega}} \right) = \infty. 
			\end{split} 
		\end{equation*} 
		
		\begin{proof}
			The proof follows the same arguments as in \cite[Lemma 2.1]{GnanasekaranCDFN24}.
		\end{proof}
	\end{lemma}
	
	Then, we proceed stating the mass limitation for $u$ and boundedness for $v,w$ as a ground for further analysis.
	
	\begin{lemma}[Mass limitation and boundedness] \label{lem:masslimBoundedness}
		Let the hypotheses of {\em Lemma \ref{lem:localex}} be complied. Then, the solution $(u,v,w)$ satisfies
		\begin{equation} \label{eq:boundU}
			\begin{split}
				\lVert u(\cdot,t) \rVert_{\Lspace{1}{\Omega}} \leq  M_1:= \max\left\{\lVert u_0 \rVert_{\Lspace{1}{\Omega}}, \frac{2 \mu_1\ \measDom  k^{\frac{k}{k+1}}}{\mu_2 (k+1)} \right\}, \quad \forall\ t\in(0, \tmax) 
			\end{split} 
		\end{equation} 
		and
		\begin{align}
			\lVert w(\cdot,t) \rVert_{\Lspace{\infty}{\Omega}} & \leq M_2:=\max\Big\{1, \|w_0\|_{\Lspace{\infty}{\Omega}}\Big\}, \hspace*{2cm}& \forall\ t\in(0, \tmax),\label{eq:boundW}\\
			\lVert v(\cdot,t) \rVert_{\Lspace{\infty}{\Omega}}&\leq M_3:=\max \Big \{\frac{\alpha}{\beta}M_2, \|v_0\|_{\Lspace{\infty}{\Omega}}\Big\}, &\forall\ t\in(0, \tmax).\label{eq:boundV}
		\end{align} 
		\begin{proof}
			By differentiating the functional $\intO u$ and using the Neumann boundary conditions, we have
			\begin{equation}\label{eq:massU}
				\frac{d}{dt}\intO u \leq \mu_1 \intO u^k - \mu_2 \intO u^{k+1} \quad \text{for all $t\in(0,\tmax)$}.
			\end{equation}
			The Young inequality gives
			\begin{equation*}
				\frac{d}{dt}\intO u \leq -\frac{\mu_2}{2} \intO u^{k+1} + \Cl[const]{c:massY}, \qquad \forall \ t\in(0,\tmax),
			\end{equation*}
			where $\Cr{c:massY}$ as  in \eqref{eqA:cmassY}.
			Moreover, the H\"older inequality leads to 
			\begin{equation*}
				\frac{d}{dt}\intO u \leq -\frac{\mu_2}{2\measDom^k} \left(\intO u\right)^{k+1} + \Cr{c:massY} \quad \text{for all}\ t\in(0,\tmax).
			\end{equation*}
			The last differential inequality, by denoting $y(t) = \intO u$ with initial condition $y(0) = \intO u_0$, provides the bound in \eqref{eq:boundU}.
			
			Further \eqref{eq:boundV} and \eqref{eq:boundW} are consequences of the parabolic maximum principle (we refer to \cite[$\S$6.4 and $\S$7.1.4]{Evans10PDEs}). 
			Indeed, the first equation entails that $u\equiv0$ is a subsolution, therefore we can assure that any solution $u$ is a nonnegative function. As a consequence, the third equation implies that $\overline{w}\equiv1$ is a supersolution, so that $0\leq w\leq1$ for all $(x,t)\in\Bar{\Omega}\times(0,\tmax)$. We exploit these information to observe that $v\equiv\frac \alpha \beta M_2$ is a supersolution for the second equation. What claimed follows by considering as well the initial conditions for $v$ and $w$.
		\end{proof}
	\end{lemma}
	
	From now on, $(u,v,w)$ indicates the local solution to problem \eqref{eq:model} provided by \emph{Lemma} \ref{lem:localex} and fulfilling the properties in {\em Lemma \ref{lem:masslimBoundedness}}.
	
	In this model \eqref{eq:model} we not only obtain the mass limitation for $u$, but further properties as inclusion of $\nabla{v}$ and $\nabla{w}$ in $\Lspace{\infty}{(0,\tmax);\Lspace{2}{\Omega}}$.
	
	\begin{lemma}[Gradient limitation] \label{lem:gradient}
		Let $k\geq 1$, $\delta > \mu_3$, and $\alpha < 2 \min\left\{\beta+\gamma,\delta - \mu_3\right\}$. Then, there exists a constant $C>0$ such that
		\begin{equation}\label{eq:boundednessgradv}
			\intO \lvert \nabla v(\cdot,t)\rvert^2 \leq C
			\quad \text{and} \quad
			\intO \lvert \nabla w(\cdot,t)\rvert^2 \leq C, \quad 
			\forall\ t\in(0,\tmax).
		\end{equation}
		
		\begin{proof}
			The thesis is derived from a differential inequality by studying the evolution of the functional
			\begin{equation} \label{eq:gradFunc}
				\varphi(t) := \tfrac{\gamma^2 M_3^2 + \delta^2 M_2^2}{\mu_2} \intO u + \intO \ngrad{v}^2 + \intO \ngrad{w}^2, \qquad \forall\ t \in (0,\tmax).
			\end{equation}
			
			For the second term in the sum, integrating by parts and by the Neumann boundary conditions, we obtain for all $t \in (0,\tmax)$ 
			\begin{equation} \label{eq:derGradv}
				\begin{split}
					\frac{d}{dt} \intO \ngrad{v}^2 
					=&\ 2 \intO \nabla v \cdot \nabla v_t = 2 \intO \nabla v \cdot \nabla (\Delta v + \alpha w - \beta v - \gamma uv) \\ 
					=& -2\intO \lvert\Delta v\rvert^2 + 2\alpha \intO \nabla v \cdot \nabla w - 2\beta \intO \ngrad{v}^2  + 2\gamma \intO uv \Delta v \\
					=& -2\intO \lvert\Delta v\rvert^2 + 2\alpha \intO \nabla v \cdot \nabla w - 2\beta \intO \ngrad{v}^2 +2 \gamma\intO (u-1)v\Delta v - 2\gamma \intO \ngrad{v}^2 
				\end{split}
			\end{equation}
			After neglecting some nonpositive terms, by considering the Cauchy-Schwarz inequality and the upper bound in \eqref{eq:boundV}, it follows
			\begin{equation}\label{eq:boundGradv}
				\begin{split}		
					\frac{d}{dt} \intO \ngrad{v}^2 \leq & -2\intO \lvert\Delta v\rvert^2 + 2\alpha \intO \nabla v \cdot \nabla w - 2(\beta +\gamma)\intO \ngrad{v}^2 + \gamma^2 \intO (u-1)^2v^2 + \intO \lvert\Delta v\rvert^2 \\
					\leq& \ 2\alpha \intO \nabla v \cdot \nabla w - 2(\beta +\gamma)\intO \ngrad{v}^2 + \gamma^2 \intO (u-1)^2v^2 \\
					\leq&\  (\alpha - 2\beta - 2\gamma) \intO \ngrad{v}^2 + \alpha \intO \ngrad{w}^2 + \gamma^2 M_3^2 \intO (u-1)^2
					\quad \text{on $(0,\tmax)$}.
				\end{split}
			\end{equation} 
			On the other hand, with an analogous reasoning as for $\intO uv\Delta v$ in \eqref{eq:derGradv} we have for all $t \in (0,\tmax)$ that
			\begin{equation} \label{eq:boundGradw}
				\begin{split}
					\frac{d}{dt} \intO \ngrad{w}^2 
					=&\ 2 \intO \nabla w \cdot \nabla w_t = 2 \intO \nabla w \cdot \nabla \left(\Delta w-\delta u w+ \mu_3 w(1-w)\right) \\
					=& -2 \intO \lvert \Delta w \rvert^2 + 2\delta \intO uw \Delta w + 2\mu_3 \intO \ngrad{w}^2 - 4\mu_3 \intO w \ngrad{w}^2 \\
					\leq&\ \delta^2 \intO w^2(u-1)^2 - 2(\delta - \mu_3) \intO \ngrad{w}^2 - 4\mu_3 \intO w \ngrad{w}^2  \\
					\leq&\ \delta^2 M_2^2 \intO (u-1)^2 - 2(\delta - \mu_3) \intO \ngrad{w}^2,
				\end{split}
			\end{equation}
			after neglecting some nonpositive terms, also in this occasion.
			
			Let us differentiate $\varphi(t)$ in \eqref{eq:gradFunc} and by considering \eqref{eq:massU}, \eqref{eq:boundGradv}, and \eqref{eq:boundGradw} we obtain
			\begin{equation} \label{eq:derPhi}
				\begin{split}
					\frac{d \varphi(t)}{dt} &\leq\  \frac{\mu_1 \Cl[const]{c:normsum} }{\mu_2} \intO u^k - \Cr{c:normsum} \intO u^{k+1} + 
					(\alpha - 2\beta - 2\gamma) \intO \ngrad{v}^2  + (\alpha - 2\delta + 2\mu_3) \intO \ngrad{w}^2  + \Cr{c:normsum} \intO (u-1)^2 \\
					&\leq\  \frac{\mu_1 \Cr{c:normsum} }{\mu_2} \intO u^k - \Cr{c:normsum} \intO u^{k+1} + 
					(\alpha - 2\beta - 2\gamma) \intO \ngrad{v}^2  + (\alpha - 2\delta + 2\mu_3) \intO \ngrad{w}^2  + \Cr{c:normsum} \intO (u+1)^2
				\end{split}
			\end{equation}
			with $\Cr{c:normsum} := \gamma^2 M_3^2 + \delta^2 M_2^2$ and $t \in (0,\tmax)$.
			
			From here we study \eqref{eq:derPhi} under the two cases: $k>1$ and $k=1$.
			
			\paragraph*{\em Case $k>1$)\ } 
			Let us denote $\Cl[const]{c:linODI}:= \max\left\{2\beta + 2\gamma -\alpha, 2\delta - 2\mu_3 -\alpha\right\} >0$, since from hypotheses $2\beta + 2\gamma -\alpha>0$ and $2\delta - 2\mu_3 -\alpha>0$. From the Young inequality for the terms $\intO u^k$ and $\intO u^2$, after adding and subtracting $\tfrac{\Cr{c:normsum}\Cr{c:linODI}}{\mu_2} \intO u$, we conclude that
			\begin{equation*} 
				\begin{split}
					\frac{d \varphi(t)}{dt} \leq&\ \Cl[const]{c:indterm}  -\tfrac{\Cr{c:normsum}\Cr{c:linODI}}{\mu_2} \intO u - (2\beta + 2\gamma -\alpha) \intO \ngrad{v}^2  - (2\delta - 2\mu_3 -\alpha) \intO \ngrad{w}^2 \\
					\leq&\  \Cr{c:indterm} - \Cr{c:linODI}\ \varphi(t), \qquad \forall\ t \in(0,\tmax),
				\end{split}
			\end{equation*}
			with $\Cr{c:indterm}$ as in \eqref{cons:indterm}. 
			
			Thus, for $C_1 := \max\left\{ \varphi(0), \frac{\Cr{c:indterm}}{\Cr{c:linODI}} \right\}$ it holds
			\begin{equation*}
				\tfrac{\gamma^2 M_3^2 + \delta^2 M_2^2}{\mu_2} \intO u + \intO \ngrad{v}^2 + \intO \ngrad{w}^2 \leq C_1, \qquad \text{on}\ (0,\tmax).
			\end{equation*}

			\paragraph*{\em Case $k=1$)\ }  In this case, we recall bound \eqref{eq:boundU}, then \eqref{eq:derPhi} becomes
			\begin{equation*}
				\begin{split}
					\frac{d \varphi(t)}{dt} \leq&\ \frac{\mu_1 \Cr{c:normsum} }{\mu_2} \intO u - \Cr{c:normsum} \intO u^2 + 
					(\alpha - 2\beta - 2\gamma) \intO \ngrad{v}^2  + (\alpha - 2\delta + 2\mu_3) \intO \ngrad{w}^2  + \Cr{c:normsum} \intO (u+1)^2 \\
					\leq&\ - \tfrac{\Cr{c:normsum}\Cr{c:linODI}}{\mu_2} \intO u - 
					(2\beta + 2\gamma - \alpha) \intO \ngrad{v}^2  - (2\delta - 2\mu_3 - \alpha) \intO \ngrad{w}^2 + \Cl[const]{c:indtermU} \\
					\leq&\ \Cr{c:indtermU} - \Cr{c:linODI} \varphi(t),
				\end{split}
			\end{equation*}
			with $\Cr{c:indtermU}$ as in \eqref{cons:indtermU}, obtaining this way a differential inequality
			valid on $(0,\tmax)$. 
			
			As in the case for $k>1$, if we take $C_2 := \max\left\{ \varphi(0), \frac{\Cr{c:indtermU}}{\Cr{c:linODI}} \right\}$ it holds
			\begin{equation*}
				\tfrac{\gamma^2 M_3^2 + \delta^2 M_2^2}{\mu_2} \intO u + \intO \ngrad{v}^2 + \intO \ngrad{w}^2 \leq C_2, \qquad \forall\ t \in(0,\tmax).
			\end{equation*}
			
			The proof is given for $C := \max \left\{C_1, C_2\right\}$.
			
		\end{proof}
	\end{lemma}
	
	The following result allows us to extend later the local-in-time existence of classical solution to globally bounded solutions. Extensibility criteria are known in the context of chemotaxis \cite{BellomoBTW15}; herein for consistency with the nomenclature we refer for details to \cite{GnanasekaranCDFN24}. 
	
	\begin{lemma}[Extensibility criterion] \label{lem:extencrit}
		Suppose there exists $p>\frac{n}{2}\geq 1$ such that
		\begin{equation*}
			\sup\limits_{t\in(0,\tmax)} \lVert u(\cdot,t) \rVert_{\Lspace{p}{\Omega}} < \infty.
		\end{equation*}
		Then we have, for all $q \geq 1$
		\begin{equation*} 
			\sup \limits_{t\in(0,\tmax)} \left( \lVert u(\cdot,t) \rVert_{\Lspace{\infty}{\Omega}} 
			+ \lVert v(\cdot,t) \rVert_{\Wspace{1}{q}{\Omega}} 
			+ \lVert w(\cdot,t) \rVert_{\Wspace{1}{q}{\Omega}} \right) < \infty. 	
		\end{equation*} 
		In particular $\tmax=\infty$ and $u \in L^\infty ((0,\infty),L^\infty(\Omega))$.
		
		\begin{proof}
			The proof follows the same reasoning as in \cite[Lemma 2.3]{GnanasekaranCDFN24}, with the main difference in the variation-of-constants formula to represent $u(\cdot,t)$ in \eqref{eq:model}. For each $t\in(0,\tmax)$ and $t_0 = \max \left\{0, t-1\right\}$, we have instead
			\begin{equation} \label{eq:NeuSGu}
				\begin{split}
					u(\cdot, t)=&\ e^{(t-t_0)\Delta} u(\cdot, t_0) - \chi \int_{t_0}^t e^{(t-s)\Delta} \nabla \cdot \left(u(\cdot,s) \nabla v(\cdot,s)\right) \ds +\mu_1 \int_{t_0}^t e^{(t-s)\Delta} u^k \left( 1 - \tfrac{\mu_2}{\mu_1} u\right) \ds. 
				\end{split}
			\end{equation}
			
			The function $z^k\left(1- \tfrac{\mu_2}{\mu_1} z\right)$, for $z\in\RR^+$ attains the maximum value $\Cl[const]{c:maxu} := \frac{k^k \mu_1^k}{(k+1)^{k+1} \mu_2^k}$ at $z = \frac{k \mu_1}{(k+1) \mu_2} > 0$. Therefore, $u^k\left(1-\tfrac{\mu_2}{\mu_1} u\right) \leq \Cr{c:maxu}$ and aside the positivity of the heat semigroup it follows that
			\begin{equation*} 
				\int_{t_0}^t e^{(t-s)\Delta} u^k \left(1-\tfrac{\mu_2}{\mu_1} u\right) \ds \leq \Cr{c:maxu}.
			\end{equation*} 
			The $\Lspace{\infty}{\Omega}$ norm of each side in \eqref{eq:NeuSGu} leads to
			\begin{equation} \label{eq:normNeuSGu}
				\lVert u(\cdot,t) \rVert_{\Lspace{\infty}{\Omega}} 
				\leq \left\lVert e^{(t-t_0)\Delta} u(\cdot,t_0) \right\rVert_{\Lspace{\infty}{\Omega}} 
				+ \chi \int_{t_0}^t \left\lVert e^{(t-s)\Delta} \nabla \cdot \left(u(\cdot,s) \nabla v(\cdot,s)\right) \right\rVert_{\Lspace{\infty}{\Omega}} \ds + \mu_1\Cr{c:maxu}.
			\end{equation}
			The bounds in \cite{GnanasekaranCDFN24} for the two first integrals in the right hand side in \eqref{eq:normNeuSGu} prove the boundedness for the left hand side.
			
			The proof follows as in \cite{GnanasekaranCDFN24}.
		\end{proof}
	\end{lemma}

	\section{A priori inequalities}
	\label{sec:aprioriest}
	

	\medskip

	With the aim to proof the boundedness of the solution $(u,v,w)$ of \eqref{eq:model} in the space $\Lspace{p}{\Omega}$, for some $p>1$, let us consider the energy functional 
	\begin{equation} \label{eq:Phidef}
		\Phi(t) = \intO u^p + \intO \ngrad{v}^{2p} + \intO \ngrad{w}^{2p}, \qquad \text{for all }\ t\in(0,\tmax),
	\end{equation} 
	with $p\in(1,\infty)$. The following lemmas pursuit that goal, treating each term in \eqref{eq:Phidef} separately before obtaining an ordinary differential equation for $\Phi(t)$.
	
	\begin{lemma} \label{lem:dFu}
		Let the hypotheses in {\em Lemma \ref{lem:gradient}} be satisfied, then it holds, for all $p>1$ and $t\in (0,\tmax)$, that
		\begin{equation} \label{eq:LemdFu}
			\dert \intO u^p + \frac{2(p-1)}{p} \intO \ngrad{u^{\frac{p}{2}}}^2 
			\leq  \frac{\chi^2 p(p-1)}{2} \intO u^p \ngrad{v}^2  - \frac{\mu_1 p}2 \intO u^{p+k} + \Cl[const]{c:dFu},
		\end{equation}
		with $\Cr{c:dFu}$ a positive computable constant.
		
		\begin{proof}
			Derivating the functional $\intO u^p$ and considering the first equation in \eqref{eq:model} we obtain on $(0,\tmax)$ that
			\begin{equation*} 
				\begin{split}
					\dert \intO u^p &= p \intO u^{p-1} u_t = p \intO u^{p-1} \left(\Delta u-\chi \nabla\cdot(u \nabla v) + \mu_1 u^k -\mu_2 u^{k+1}\right) \\
					&= -p(p-1)\intO u^{p-2} \ngrad{u}^{2} + \chi p(p-1) \intO u^{p-1}\nabla u \cdot \nabla v + \mu_1 p \intO u^{p+k-1} - \mu_2 p \intO u^{p+k} \\
					&\leq -\frac{p(p-1)}{2}\intO u^{p-2} \ngrad{u}^{2} + \frac{\chi^2 p(p-1)}{2} \intO u^{p} \lvert \nabla v \rvert^2 + \mu_1 p \intO u^{p+k-1} - \mu_2 p \intO u^{p+k} \\
					&\leq -\frac{2(p-1)}{p} \intO \ngrad{u^{\frac{p}{2}}}^2 + \frac{\chi^2 p(p-1)}{2} \intO u^p \ngrad{v}^2 - \frac{\mu_2 p}2 \intO u^{p+k} + \Cr{c:dFu},
				\end{split}
			\end{equation*}
			after applying the Young inequality with $\Cr{c:dFu}$ defined in \eqref{eqA:dFu}.
		\end{proof}
	\end{lemma}
	
	\begin{lemma} \label{lem:dFgradv}
		Let the hypotheses in {\em Lemma \ref{lem:gradient}} be satisfied, then for any $p>1$, certain $\Cl[eps]{e:dFgradvI}, \Cl[eps]{e:dFgradvII} > 0$, $\eta \in (0,p-1)$, and computable constants $\Cl[const]{c:tempI}(n,\eta), \Cl[const]{c:tempII}(\Cr{e:dFgradvI},\Cr{e:dFgradvII}) >0$, it holds, for all $t \in (0,\tmax)$, that
		\begin{multline} \label{eq:LemdFgradv}
			\dert \intO \ngrad{v}^{2p} \leq - p \intO \ngrad{v}^{2p-2} \lvert D^2 v\rvert^2 + \Cr{e:dFgradvI} \intO \ngrad{v}^{2(p+1)} + \Cr{e:dFgradvII} \intO \ngrad{w}^{2(p+1)} \\ + \Cr{c:tempI}(n,\eta) \gamma^2 M_3^2 \intO u^2 \ngrad{v}^{2p-2} + \Cr{c:tempII}(\Cr{e:dFgradvI},\Cr{e:dFgradvII}).
		\end{multline} 
		
		\begin{proof}
			Similarly to what done before, after substituting $v_t$ from \eqref{eq:model} and integrating by parts some integral terms on $(0,\tmax)$, we have 
			\begin{equation*}
				\begin{split}
					\dert \intO \ngrad{v}^{2p} =&\ p \intO \ngrad{v}^{2p-2} 2\nabla v \cdot \nabla v_t 
					= p \intO \ngrad{v}^{2p-2} 2\nabla v \cdot \nabla \left(\Delta v + \alpha w - \beta v - \gamma uv\right) \\
					=&\ p\intO \ngrad{v}^{2p-2} \left(\Delta \ngrad{v}^2 - 2\lvert D^2 v\rvert^2 +2\alpha\nabla v \cdot \nabla w - 2\beta \ngrad{v}^2 - 2\gamma \nabla v \cdot \nabla (uv) \right) \\
					=&\ p\int_{\partial\Omega} \ngrad{v}^{2p-2} \ \partialn{\ngrad{v}^2} -p(p-1) \intO \ngrad{v}^{2p-4} \ngrad{\ngrad{v}^2}^2 - 2p \intO \ngrad{v}^{2p-2} \lvert D^2 v\rvert^2 \\
					&\  + 2\alpha p \intO \ngrad{v}^{2p-2} \nabla v \cdot \nabla w - 2\beta p \intO \ngrad{v}^{2p} - 2\gamma p \intO \ngrad{v}^{2p-2} \nabla v \cdot \nabla (uv) \\
					=&\ p\int_{\partial\Omega} \ngrad{v}^{2p-2} \ \partialn{\ngrad{v}^2} -p(p-1) \intO \ngrad{v}^{2p-4} \ngrad{\ngrad{v}^2}^2 - 2p \intO \ngrad{v}^{2p-2} \lvert D^2 v\rvert^2 \\
					&\  + 2\alpha p \intO \ngrad{v}^{2p-2} \nabla v \cdot \nabla w - 2\beta p \intO \ngrad{v}^{2p} + 2\gamma p \intO uv \ngrad{v}^{2p-2} \Delta v \\
					& + 2\gamma p(p-1) \intO uv \ngrad{v}^{2p-4} \nabla v \cdot \nabla \ngrad{v}^2,
				\end{split}
			\end{equation*}
			with the second line justified by the identity \eqref{eq:idLapNormGrad}.
			
			After ignoring the nonpositive fifth term, we use the inequality in \eqref{ineq:boundaryGradIneq} for some $\eta \in (0,p-1)$ to treat the boundary integral. Successively, bounds in \eqref{eq:boundV} and \eqref{eq:boundednessgradv}, aside the Cauchy-Schwarz inequality lead then to
			\begin{equation*} 
				\begin{split}
					\dert \intO \ngrad{v}^{2p} \leq&\ -p\left(p-1-\eta\right)\intO \ngrad{v}^{2p-4}\ngrad{\ngrad{v}^2}^2 + p C_\eta \left(\intO \ngrad{v}^2\right)^p - 2p \intO \ngrad{v}^{2p-2} \lvert D^2 v\rvert^2 \\
					&  + 2\alpha p \intO \ngrad{v}^{2p-2} \nabla v \cdot \nabla w  + 2\gamma p M_3 \intO u \ngrad{v}^{2p-2} \Delta v + 2\gamma p(p-1) M_3 \intO u \ngrad{v}^{2p-4} \nabla v \cdot \nabla \ngrad{v}^2 \\
					\leq&\ -p\left(p-1-\eta\right)\intO \ngrad{v}^{2p-4}\ngrad{\ngrad{v}^2}^2 - 2p \intO \ngrad{v}^{2p-2} \lvert D^2 v\rvert^2 + 2\alpha p \intO \ngrad{v}^{2p-1} \lvert\nabla w\rvert \\
					&   + 2\gamma p M_3 \intO u \ngrad{v}^{2p-2} \lvert\Delta v\rvert + 2\gamma p(p-1) M_3 \intO u \ngrad{v}^{2p-3}  \ngrad{\ngrad{v}^2} + \Cl[const]{c:LankeitineqU}(\eta)
				\end{split}
			\end{equation*} 
			for all $t \in (0,\tmax)$, where $\Cr{c:LankeitineqU}(\eta) := p\ C_\eta C^p$, for the constants $C$ obtained in \eqref{eq:boundednessgradv} and $C_\eta$ in \eqref{ineq:boundaryGradIneq}.
			
			After applying the Young inequality for 
			\begin{math}
				u \ngrad{v}^{2p-2} \lvert\Delta v\rvert 
				%
			\end{math}
			and 
			\begin{math}
				u \ngrad{v}^{2p-3}  \ngrad{\ngrad{v}^2} 
			\end{math},
			it follows from \eqref{ineq:LapHessian} and \eqref{ineq:gradVgradW} that
			\begin{equation*} 
				\begin{split}
					\dert \intO \ngrad{v}^{2p} \leq&\ - p \intO \ngrad{v}^{2p-2} \lvert D^2 v\rvert^2 + 2\alpha p \intO \ngrad{v}^{2p-1} \lvert\nabla w\rvert + \Cr{c:tempI}(n,\eta) \gamma^2 M_3^2 \intO u^2 \ngrad{v}^{2p-2} + \Cr{c:LankeitineqU}(\eta) \\
					\leq& - p \intO \ngrad{v}^{2p-2} \lvert D^2 v\rvert^2 
					+ \Cr{e:dFgradvI} \intO \ngrad{v}^{2(p+1)} 
					+ \Cr{e:dFgradvII} \intO \ngrad{w}^{2(p+1)}  \\ 
					& + \Cr{c:tempI}(n,\eta) \gamma^2 M_3^2 \intO u^2 \ngrad{v}^{2p-2} + \Cr{c:tempII}(\Cr{e:dFgradvI},\Cr{e:dFgradvII},\eta) 
				\end{split}
			\end{equation*}
			with $\Cr{c:tempI}(n,\eta) := p \left(n + \frac{(p-1)^2}{p-1-\eta}\right)$ and $\Cr{c:tempII}(\Cr{e:dFgradvI},\Cr{e:dFgradvII},\eta)$ defined in \eqref{eqA:dFgradv}.
			
		\end{proof}
	\end{lemma}
	
	\begin{lemma} \label{lem:dFgradw}
		Let the hypotheses in {\em Lemma \ref{lem:gradient}} be satisfied, then for $p>1$, a computable constant $\Cl[const]{c:tempIII} >0$ (with $\Cr{c:tempI}(n,\eta)$ defined as in \eqref{eq:LemdFgradv}, and $\eta \in (0,p-1)$), it holds, for all $t \in (0,\tmax)$, that
		\begin{equation} \label{eq:LemdFgradw}
			\dert \intO \ngrad{w}^{2p} \leq - p \intO \ngrad{w}^{2p-2} \lvert D^2 w\rvert^2 + \Cl[eps]{e:ineqGradW} \intO \ngrad{\ngrad{w}^{p}}^2 + \Cr{c:tempI}(n,\eta) \delta^2 M_2^2 \intO u^2 \ngrad{w}^{2p-2}  + \Cr{c:tempIII}(\Cr{e:ineqGradW}).
		\end{equation}
		
		\begin{proof}
			Following an analogous procedure as before, for the third integral term in \eqref{eq:Phidef}, we have for all $t\in(0,\tmax)$
			\begin{equation*} 
				\begin{split}
					\dert \intO \ngrad{w}^{2p} =&\ p \intO \ngrad{w}^{2p-2} 2\nabla w \cdot \nabla \left(\Delta w - \delta u w +\mu_3 (1 - w)\right) \\
					=&\ p \intO \ngrad{w}^{2p-2} \left(\Delta \ngrad{w}^2 - 2\lvert D^2 w\rvert^2\right) -2p\delta \intO \ngrad{w}^{2p-2} \nabla w \cdot \nabla (uw) \\
					& + 2p\mu_3 \intO \ngrad{w}^{2p} - 4p\mu_3 \intO w \ngrad{w}^{2p} \\
					\leq& - p \intO \ngrad{w}^{2p-2} \lvert D^2 w\rvert^2 + 2p\mu_3 \intO \ngrad{w}^{2p} + \Cr{c:tempI}(n,\eta) \delta^2 M_2^2 \intO u^2 \ngrad{w}^{2p-2}  + \Cr{c:LankeitineqU}(\eta) \\
					\leq& - p \intO \ngrad{w}^{2p-2} \lvert D^2 w\rvert^2 + \frac{1}{p^2} \intO \ngrad{\ngrad{w}^{p}}^2 + \Cr{c:tempI}(n,\eta) \delta^2 M_2^2 \intO u^2 \ngrad{w}^{2p-2}  + \Cr{c:tempIII}(\eta),
				\end{split}
			\end{equation*}
			after using \eqref{ineq:gradW} and with 
			$\Cr{c:tempIII}(\Cr{e:ineqGradW},\eta)$ as defined in \eqref{eqA:dFgradw}.
			
		\end{proof}
	\end{lemma}
	
	Finally, let us use the  estimates obtained from the previous three lemmas to prove the boundedness of the functional defined in \eqref{eq:Phidef}.
	
	\begin{lemma} \label{lem:boundPhi}
		Let the hypotheses of {\em Lemma \ref{lem:gradient}} be satisfied, $p>1$ and either $k>1$ or $k=1$ with
		\begin{equation}  \label{eq:hypoMup}
			\mu_2 > \ A_1(p,n) \chi^{2+\frac2p} M_3^\frac2p
			+ A_2(p,n) \gamma^{p+1}  M_3^{2p} + A_3(p,n) \delta^{p+1} M_2^{2p}
		\end{equation}
		hold, where 
		\begin{align}
			A_1(p,n) &:=\ \tfrac{(p-1) p^2}{2(p+1)} \left(\tfrac{6\ (p-1) \left(4 p^2 + n\right)}{p+1}\right)^{1/p} \label{eq:defA1} \\
			A_2(p,n) &:=\ \tfrac{2\ \Cr{c:tempI}(n,0)}{p+1} \left(\tfrac{12\ \Cr{c:tempI}(n,0) \left(4 p^2 + n\right)}{p (p+1)}\right)^{\frac{p-1}{2}} \label{eq:defA2} \\
			A_3(p,n) &:=\ \tfrac{2\ \Cr{c:tempI}(n,0)}{p+1} \left(\tfrac{16\ \Cr{c:tempI}(n,0) (p-1) \left(4 p^2 + n\right)}{p (p+1)}\right)^{\frac{p-1}{2}} . \label{eq:defA3}
		\end{align}
		(with $\Cr{c:tempI}(n,\eta)$ defined as in \eqref{eq:LemdFgradv}).
		
		Then, there exists a constant $\widetilde{C}$ such that 
		\begin{equation*}
			\intO u^p + \intO \ngrad{v}^{2p} + \intO \ngrad{w}^{2p} \leq \widetilde{C} \qquad \forall \ t \in (0,\tmax).
		\end{equation*}
		
		\begin{proof}	
			We start by combining \eqref{eq:LemdFu}, \eqref{eq:LemdFgradv}, and \eqref{eq:LemdFgradw} to obtain an estimate for $\dert \Phi(t)$ on $(0,\tmax)$ from \eqref{eq:Phidef}
			\begin{equation} \label{eq:dPhibeforeK}
				\begin{split}
					\dert \Phi(t) =&\ \dert \intO u^p + \dert \intO \ngrad{v}^{2p}  + \dert \intO \ngrad{w}^{2p} \\
					\leq& -\frac{2(p-1)}{p} \intO \ngrad{u^{\frac{p}{2}}}^2  - p \intO \ngrad{v}^{2p-2} \lvert D^2 v\rvert^2 - p \intO \ngrad{w}^{2p-2} \lvert D^2 w\rvert^2 - \frac{\mu_2 p}2 \intO u^{p+k}  \\
					& + \frac{\chi^2 p(p-1)}{2} \intO u^p \ngrad{v}^2 + \Cr{c:tempI}(n,\eta) \gamma^2 M_3^2 \intO u^2 \ngrad{v}^{2p-2} + \Cr{c:tempI}(n,\eta) \delta^2 M_2^2 \intO u^2 \ngrad{w}^{2p-2}		
					\\ 
					&  + \Cr{e:dFgradvI} \intO \ngrad{v}^{2(p+1)} + \Cr{e:dFgradvII} \intO \ngrad{w}^{2(p+1)} + \Cr{e:ineqGradW} \intO \ngrad{\ngrad{w}^{p}}^2 + \Cl[const]{c:tempIV}(\Cr{e:dFgradvI},\Cr{e:dFgradvII},\eta), 
				\end{split}
			\end{equation}
			with $\Cr{c:tempIV}(\Cr{e:dFgradvI},\Cr{e:dFgradvII},\eta):= \Cr{c:dFu} + \Cr{c:tempII}(\Cr{e:dFgradvI},\Cr{e:dFgradvII},\eta) + \Cr{c:tempIII}(\Cr{e:ineqGradW},\eta)$.
			
			From here on we distinguish between two cases: $k>1$ and $k=1$.
			
			\paragraph*{\em Case $k>1$)\ } 
			
			From the inequalities \eqref{ineq:uPgradv2} and \eqref{ineq:u2gradvDPMD}, after setting $\Cr{e:ineqGradW}=\frac{1}{4p}$, we obtain for $t\in (0,\tmax)$ that
			\begin{equation} \label{eq:dPhiKgeqUi}
				\begin{split}
					\dert \Phi(t) \leq&  -\frac{2(p-1)}{p} \intO \ngrad{u^{\frac{p}{2}}}^2  - p \intO \ngrad{v}^{2p-2} \lvert D^2 v\rvert^2 - p \intO \ngrad{w}^{2p-2} \lvert D^2 w\rvert^2 \\
					& + \left(\Cl[eps]{e:uYoungI} + \Cl[eps]{e:uYoungII} + \Cl[eps]{e:uYoungIII} - \tfrac{\mu_2 p}2\right) \intO u^{p+k} 
					+ (\Cr{e:dFgradvI} + \Cl[eps]{e:gradvYoungI} + \Cl[eps]{e:gradvYoungII}) \intO \ngrad{v}^{2(p+1)} 
					+ (\Cr{e:dFgradvII} + \Cl[eps]{e:gradwYoungI}) \intO \ngrad{w}^{2(p+1)} \\
					& + \frac{1}{4p} \intO \ngrad{\ngrad{w}^{p}}^2 + \Cl[const]{c:indKgequI}(\Cr{e:dFgradvI}, \Cr{e:dFgradvII}, \Cr{e:uYoungII}, \Cr{e:uYoungIII},\Cr{e:gradvYoungII}, \Cr{e:gradwYoungI}, \eta)
				\end{split}
			\end{equation}
			where 
			\begin{math}
				\Cr{c:indKgequI}(\Cr{e:dFgradvI}, \Cr{e:dFgradvII}, \Cr{e:uYoungII}, \Cr{e:uYoungIII},\Cr{e:gradvYoungII}, \Cr{e:gradwYoungI}, \eta) := \Cr{c:tempIV}(\Cr{e:dFgradvI},\Cr{e:dFgradvII},\eta) + \Cl[const]{c:uPgradv2}(\Cr{e:uYoungIII},\Cr{e:gradwYoungI}) + \Cl[const]{c:u2gradvDPMD}(\Cr{e:uYoungII},\Cr{e:gradvYoungII},\eta) + \Cl[const]{c:u2gradwDPMD}(\Cr{e:uYoungIII},\Cr{e:gradwYoungI},\eta)
			\end{math}
			with the computable constants $\Cr{c:uPgradv2}(\Cr{e:uYoungI},\Cr{e:gradvYoungI})$, $\Cr{c:u2gradvDPMD}(\Cr{e:uYoungII},\Cr{e:gradvYoungII},\eta)$, and $\Cr{c:u2gradwDPMD}(\Cr{e:uYoungIII},\Cr{e:gradwYoungI},\eta)$ obtained from the Young inequality.
			
			
			Let us choose one among the infinite positive solutions of the linear system of equations
			\begin{equation} \label{eq:systemK}
				\begin{cases}
					\Cr{e:uYoungI} + \Cr{e:uYoungII} + \Cr{e:uYoungIII} &= \frac{\mu_2 p}2 \\
					\Cr{e:dFgradvI} + \Cr{e:gradvYoungI} + \Cr{e:gradvYoungII} &= \frac{p}{4(4p^2 + n)M_3^2} \\
					\Cr{e:dFgradvII} + \Cr{e:gradwYoungI} &= \frac{p}{8(4p^2 + n)M_2^2} 
				\end{cases}.
			\end{equation}
			Then, thanks to \eqref{ineq:gradHessian} for an arbitrary value of $\eta \in (0,p-1)$ and for all $t\in(0,\tmax)$, the inequality in \eqref{eq:dPhiKgeqUi} becomes
			\begin{equation} \label{eq:dPhiKgeqUii}
				\begin{split}
					\dert \Phi(t) \leq&  -\frac{2(p-1)}{p} \intO \ngrad{u^{\frac{p}{2}}}^2  - \frac{p}{2}  \intO \ngrad{v}^{2p-2} \lvert D^2 v\rvert^2 - \frac{3p}{4} \intO \ngrad{w}^{2p-2} \lvert D^2 w\rvert^2 \\
					& + \frac{1}{4p} \intO \ngrad{\ngrad{w}^{p}}^2 + \Cr{c:indKgequI} \\
					\leq& -\frac{2(p-1)}{p} \intO \ngrad{u^{\frac{p}{2}}}^2  - \frac{1}{2p} \intO \ngrad{\ngrad{v}^{p}}^2 - \frac{1}{2p} \intO \ngrad{\ngrad{w}^{p}}^2 + \Cr{c:indKgequI},
				\end{split}
			\end{equation}
			due to the inequality \eqref{ineq:NgradNgrad}.
			
			From the Gagliardo-Nirenberg inequality and the bounds in \eqref{eq:boundednessgradv} follow the inequalities on $(0,\tmax)$
			\begin{equation} \label{eq:GNforUp}
				\intO u^p \leq 2 C_{GN}^2 M_1^{1-\theta} \left(\intO \ngrad{u^{\frac{p}{2}}}^2\right)^{\theta} + 2 C_{GN}^2 M_1^2, 
			\end{equation}
			\begin{equation} \label{eq:GNforGradvDp}
				\intO \ngrad{v}^{2p} \leq 2C_{GN}^2 C^{1-\theta} \left(\intO \ngrad{\ngrad{v}^p}^2\right)^{\theta} + 2 C_{GN}^2 C^2,
			\end{equation}
			and analogous to \eqref{eq:GNforGradvDp} for $\intO \ngrad{w}^{2p}$, 
			where 
			$C$ is the constant in {\em Lemma \ref{lem:gradient}} and $\theta = \frac{p-1}{p-1 + \frac1n}\in(0,1)$.
			
			Therefore, recalling $(a+b)^\frac{1}{\theta} \leq 2^{\frac{1-\theta}{\theta}} \left(a^\frac{1}{\theta} + b^\frac{1}{\theta}\right)$ provided $\theta<1$ and $a,b\geq 0$, and rearranging terms for each inequality \eqref{eq:GNforUp} and \eqref{eq:GNforGradvDp} we have on $(0,\tmax)$
			\begin{equation} \label{eq:GNforUpneg}
				- \intO \ngrad{u^{\frac{p}{2}}}^2 \leq - \Cl[const]{c:GNforUp} \left(\intO u^p\right)^{\frac{1}{\theta}} + \Cl[const]{c:GNforUpIndT}
			\end{equation}
			and 
			\begin{equation} \label{eq:GNforGradvDpneg}
				- \intO \ngrad{\ngrad{v}^p}^2 \leq - \Cl[const]{c:GNforGradvDp} \left(\intO \ngrad{v}^{2p}\right)^{\frac{1}{\theta}} + \Cl[const]{c:GNforGradvDpIndT},
			\end{equation}
			with a similar expression to the latter for the function $w$, being 
			\begin{equation*}
				\begin{split}
					\Cr{c:GNforUp} &:= \frac{1}{\left( 2 C_{GN}^2 M_1^{1-\theta} \right)^{\frac{1}{\theta}}}, 
					\quad 
					\Cr{c:GNforUpIndT} := \left(2 C_{GN}^2 M_1^2\right)^{\frac{1}{\theta}} \ \Cr{c:GNforUp}, \\
					\Cr{c:GNforGradvDp} &:= \frac{1}{\left(2C_{GN}^2 C^{1-\theta}\right)^{\frac{1}{\theta}}},
					\quad 
					\Cr{c:GNforGradvDpIndT} := \left(2C_{GN}^2 C^{1-\theta}\right)^{\frac{1}{\theta}} \ \Cr{c:GNforGradvDp}.
				\end{split}
			\end{equation*}
			
			Inserting the last inequalities \eqref{eq:GNforUpneg} and \eqref{eq:GNforGradvDpneg} into \eqref{eq:dPhiKgeqUii} we conclude
			\begin{equation} \label{eq:obtODIenergy}
				\begin{split}
					\dert \Phi(t) \leq& -\frac{2(p-1)}{p} \Cr{c:GNforUp} \left(\intO u^p\right)^{\frac{1}{\theta}}
					- \frac{\Cr{c:GNforGradvDp}}{2p}  \left(\intO \ngrad{v}^{2p}\right)^{\frac{1}{\theta}}  
					- \frac{\Cr{c:GNforGradvDp}}{2p}  \left(\intO \ngrad{w}^{2p}\right)^{\frac{1}{\theta}} \\ 
					&+ \frac{2(p-1)}{p} \Cr{c:GNforUpIndT}  + \frac{\Cr{c:GNforGradvDpIndT}}{p} + \Cr{c:indKgequI} \\
					\leq&  -\Cl[const]{c:coefPhi} \left( \left(\intO u^p\right)^{\frac{1}{\theta}} + \left(\intO \ngrad{v}^{2p}\right)^{\frac{1}{\theta}} + \left(\intO \ngrad{w}^{2p}\right)^{\frac{1}{\theta}} \right) + \Cl[const]{c:indPhi} \\
					\leq& -\frac{\Cr{c:coefPhi}}{3^{\frac{1-\theta}{\theta}}} \Phi(t)^{\frac{1}{\theta}} + \Cr{c:indPhi}, \qquad \forall \ t\in (0,\tmax),
				\end{split}
			\end{equation}
			where $\Cr{c:coefPhi} := \min \left\{\frac{2(p-1)\Cr{c:GNforUp}}{p}, \frac{\Cr{c:GNforGradvDp}}{2p}\right\}$ and $\Cr{c:indPhi}:=\frac{2(p-1)\Cr{c:GNforUpIndT}}{p} + \frac{\Cr{c:GNforGradvDpIndT}}{p} + \Cr{c:indKgequI}$, and the last inequality, justified by the generalized mean inequality $(a+b+c)^\frac{1}{\theta} \leq 3^{\frac{1-\theta}{\theta}} \left(a^\frac{1}{\theta} + b^\frac{1}{\theta} + c^\frac{1}{\theta}\right)$, with $\theta<1$ and $a,b,c\geq 0$. Thus, we obtain the differential inequality 
			\begin{equation} \label{eq:ineqDifPhi}
				\begin{dcases}
					\Phi'(t) \leq -\frac{\Cr{c:coefPhi}}{{3^{\frac{1-\theta}{\theta}}}} \Phi(t)^{\frac{1}{\theta}} + \Cr{c:indPhi}, & \qquad \forall \ t\in(0,\tmax) \\
					\Phi(0) = \intO u_0^{2p} + \intO \ngrad{v_0}^{2p} + \intO \ngrad{w_0}^{2p} .
				\end{dcases}
			\end{equation}
			
			From the comparison principle, it follows the upper bound
			\begin{equation} \label{eq:boundPhi}
				\intO u^p + \intO \ngrad{v}^{2p} + \intO \ngrad{w}^{2p} \leq C_3 := \max \left\{\Phi(0), 3^{\frac{1}{1 + n (p-1)}} \ \left( \tfrac{\Cr{c:indPhi}}{\Cr{c:coefPhi}}\right)^{\frac{p-1}{p-1 + \frac1n}}\right\}, \qquad \forall \ t\in (0,\tmax).
			\end{equation}
			
			
			\paragraph*{\em Case $k=1$)\ } 
			
			Back to \eqref{eq:dPhibeforeK}, this time the inequalities \eqref{ineq:uPgradv2KO} and \eqref{ineq:u2gradvDPMDKO}, after choosing $\Cr{e:ineqGradW}=\frac{1}{4p}$, provide
			\begin{equation} \label{eq:derPhiK1}
				\begin{split}
					%
					\dert \Phi(t) \leq&  -\frac{2(p-1)}{p} \intO \ngrad{u^{\frac{p}{2}}}^2  - p \intO \ngrad{v}^{2p-2} \lvert D^2 v\rvert^2 - p \intO \ngrad{w}^{2p-2} \lvert D^2 w\rvert^2 \\
					& + (\Cl[eps]{e:uYoungI1} + \Cl[eps]{e:uYoungII1} + \Cr{e:dFgradvI}) \intO \ngrad{v}^{2(p+1)} + (\Cl[eps]{e:uYoungIII1} + \Cr{e:dFgradvII}) \intO \ngrad{w}^{2(p+1)} + \frac{1}{4p^2} \intO \ngrad{\ngrad{w}^{p}}^2 \\ 
					& + \left(\Cl[const]{c:gradvYoungI1}(\Cr{e:uYoungI1}) + \Cl[const]{c:gradvYoungII1}(\Cr{e:uYoungII1},n,\eta) + \Cl[const]{c:gradwYoungI1}(\Cr{e:uYoungIII1},n,\eta) - \tfrac{\mu_2 p}2\right) \intO u^{p+1} + \Cr{c:tempIV}(\Cr{e:dFgradvI}, \Cr{e:dFgradvII}, \eta), \quad \forall \ t\in (0,\tmax), 
				\end{split}
			\end{equation}
			where the computable constants $\Cr{c:gradvYoungI1}(\Cr{e:uYoungI1}), \Cr{c:gradvYoungII1}(\Cr{e:uYoungII1},n,\eta)$, and $\Cr{c:gradwYoungI1}(\Cr{e:uYoungIII1},n,\eta)$ are obtained from the Young inequa\-lity being
			\begin{equation} \label{eq:epsYoungI1}
				\Cr{c:gradvYoungI1}(\Cr{e:uYoungI1}) := \left(\tfrac{\chi^2 p(p-1)}{2}\right)^{\tfrac{p+1}{p}}\ \tfrac{p}{(p+1)\left(\Cr{e:uYoungI1} (p+1) \right)^{\frac{1}{p}}} ,
			\end{equation}
			\begin{equation} \label{eq:epsYoungII1}
				\Cr{c:gradvYoungII1}(\Cr{e:uYoungII1},n,\eta) := \left(\Cr{c:tempI}(n,\eta) \gamma^2 M_3^2\right)^{\tfrac{p+1}{2}} \ \tfrac{2}{p+1} \left( \tfrac{\Cr{e:uYoungII1}(p+1)}{p-1}\right)^{\frac{1-p}{2}} ,
			\end{equation}
			and
			\begin{equation} \label{eq:epsYoungIII1}
				\Cr{c:gradwYoungI1}(\Cr{e:uYoungIII1},n,\eta) := \left(\Cr{c:tempI}(n,\eta) \delta^2 M_2^2\right)^{\tfrac{p+1}{2}} \ \tfrac{2}{p+1} \left( \tfrac{\Cr{e:uYoungIII1}(p+1)}{p-1}\right)^{\frac{1-p}{2}}.
			\end{equation}
			
			Let us choose analogously one among the infinite positive solutions of the linear system of equations
			\begin{equation} \label{eq:systemK1}
				\begin{cases}
					\Cr{e:uYoungI1} + \Cr{e:uYoungII1} + \Cr{e:dFgradvI} &= \frac{p}{4(4p^2+n)M_3^2} \\
					\Cr{e:uYoungIII1} + \Cr{e:dFgradvII} &= \frac{p}{8(4p^2 + n)M_2^2}
				\end{cases}.
			\end{equation} 
			In addition, we require that
			\begin{math}
				\Cr{c:gradvYoungI1}(\Cr{e:uYoungI1}) + \Cr{c:gradvYoungII1}(\Cr{e:uYoungII1},n,\eta) + \Cr{c:gradwYoungI1}(\Cr{e:uYoungIII1},n,\eta) \leq \frac{\mu_2 p}2.
			\end{math}
			That inequality is valid whenever the parameter $\mu_2$ is chosen large enough once a solution of the system of equations \eqref{eq:systemK1} is selected.
			
			Let us consider the particular solution satisfying $\Cr{e:uYoungI1} = \Cr{e:uYoungII1} = \Cr{e:dFgradvI}$ and $\Cr{e:uYoungIII1} = \Cr{e:dFgradvII}$. Evaluating in \eqref{eq:epsYoungI1}, \eqref{eq:epsYoungII1}, and \eqref{eq:epsYoungIII1}, we obtain the restriction on $\mu_2$ reading $\mu_2 \geq f(p,\eta)$ with
			\begin{multline*} 
				f(p,\eta) := \ \tfrac{(p-1) p^2}{2(p+1)} \left(\tfrac{6\ (p-1) \left(4 p^2 + n\right)}{p+1}\right)^{1/p} \chi^{2+\frac2p} M_3^\frac2p
				+ \tfrac{2\ \Cr{c:tempI}(n,\eta)}{p+1} \left(\tfrac{12\ \Cr{c:tempI}(n,\eta) \left(4 p^2 + n\right)}{p (p+1)}\right)^{\frac{p-1}{2}} \gamma^{p+1} M_3^{2p} \\ 
				+ \tfrac{2\ \Cr{c:tempI}(n,\eta)}{p+1} \left(\tfrac{16\ \Cr{c:tempI}(n,\eta) (p-1) \left(4 p^2 + n\right)}{p (p+1)}\right)^{\frac{p-1}{2}} \delta^{p+1} M_2^{2p}.
			\end{multline*}
			From the definition of $f(p,\eta)$, we can rewrite the hypothesis \eqref{eq:hypoMup} as
			\begin{equation} \label{eq:restrictionMup}
				\mu_2 - f(p,0) = \mu_2 - A_1(p,n) \chi^{2+\frac2p} M_3^\frac2p
				- A_2(p,n) \gamma^{p+1}  M_3^{2p} - A_3(p,n) \delta^{p+1} M_2^{2p} > 0,
			\end{equation}
			with $A_1(p,n), A_2(p,n)$, and $A_3(p,n)$ already defined in \eqref{eq:defA1}, \eqref{eq:defA2}, and \eqref{eq:defA3} respectively. Being $\mu_2 - f(p,\eta)$ a continuous function on both variables $p$ and $\eta$, then there exists an $\eta \in (0,p-1)$, such that $\mu_2 - f(p,\eta) \geq 0$. That implies that the coefficient of the term $\intO u^{p+1}$ in \eqref{eq:derPhiK1} is nonpositive. Therefore, after neglecting that term and with analogous reasoning as from \eqref{eq:dPhiKgeqUii} to \eqref{eq:obtODIenergy}, the inequalities in \eqref{eq:GNforUpneg} and \eqref{eq:GNforGradvDpneg} imply that
			\begin{equation*}
				\begin{split}
					\dert \Phi(t) \leq& -\frac{2(p-1)}{p} \intO \ngrad{u^{\frac{p}{2}}}^2 - \frac{1}{2p} \intO \ngrad{\ngrad{v}^{p}}^2 - \frac{1}{2p} \intO \ngrad{\ngrad{w}^{p}}^2 + \Cr{c:tempIV}(\eta) \\
					\leq& -\Cr{c:coefPhi} \left( \left(\intO u^p\right)^{\frac{1}{\theta}} + \left(\intO \ngrad{v}^{2p}\right)^{\frac{1}{\theta}} + \left(\intO \ngrad{w}^{2p}\right)^{\frac{1}{\theta}} \right) + \Cl[const]{c:indPhiKi} \\
					\leq& -\frac{\Cr{c:coefPhi}}{3^{\frac{1-\theta}{\theta}}} \Phi(t)^{\frac{1}{\theta}} + \Cr{c:indPhiKi}, \qquad \forall \ t\in (0,\tmax),
				\end{split}
			\end{equation*}
			with $\Cr{c:indPhiKi} := \frac{2(p-1)\Cr{c:GNforUpIndT}}{p} + \frac{\Cr{c:GNforGradvDpIndT}}{p} + \Cr{c:tempIV}(\eta)$. 
			
			With the same argument as in \eqref{eq:ineqDifPhi} we obtain the upper bound
			\begin{equation} \label{eq:boundPhiKi}
				\intO u^p + \intO \ngrad{v}^{2p} + \intO \ngrad{w}^{2p} \leq C_4 := \max \left\{\Phi(0), 3^{\frac{1}{1 + n (p-1)}} \ \left( \tfrac{\Cr{c:indPhiKi}}{\Cr{c:coefPhi}}\right)^{\frac{p-1}{p-1 + \frac1n}}\right\}, \qquad \forall \ t\in (0,\tmax).
			\end{equation}
			
			It is enough to consider as the upper bound in the thesis $\widetilde{C} = \max\{C_3,C_4\}$, with $C_3$ defined in \eqref{eq:boundPhi} and $C_4$ in \eqref{eq:boundPhiKi}.
			
		\end{proof}
	\end{lemma}

	
	With the previous result we are able to provide a proof to {\em Theorem \ref{theo:globalexist}}.
	
	\begin{proof}
		Let us divide the proof into two cases: $k>1$ and $k=1$, which rely on different hypotheses.
		
		For $k>1$, the hypotheses in \eqref{eq:theocond} agree with the requirements of {\em Lemma \ref{lem:gradient}}. Thus, {\em Lemma \ref{lem:boundPhi}} is valid and, as a consequence of {\em Lemma \ref{lem:extencrit}}, it follows the thesis.
		
		On the other hand, for $k=1$, let us consider the function $\mu_2 - f(p,0)$ already defined in \eqref{eq:restrictionMup}, which is continuous for the variable $p$ by the definition of $A_1(p,n), A_2(p,n)$, and $A_3(p,n)$. With abuse of notation, let us define those constants as: $A_1(n) = A_1(\frac{n}{2},n), A_2(n) = A_2(\frac{n}{2},n)$, and $A_3(n) = A_3(\frac{n}{2},n)$. The hypothesis in \eqref{eq:muRaf} states that $\mu_2 - f(\frac{n}{2},0) > 0$, then there exists $p > \frac{n}{2}$ such that the hypothesis \eqref{eq:hypoMup} of {\em Lemma \ref{lem:boundPhi}} holds. 		
		Then, {\em Lemma \ref{lem:gradient}} follows from the hypotheses \eqref{eq:theocond} and consequently {\em Lemma \ref{lem:boundPhi}} applies again. 
	\end{proof}
	
	\subsection{Further reflections}
	
	Let us consider the case $\alpha=\beta=0$, then the studied model becomes
	\begin{equation} \label{eq:modelWithoutAB}
		\begin{dcases}
			u_t = \Delta u-\chi \nabla\cdot(u \nabla v) + \mu_1 u^k - \mu_2 u^{k+1},\hspace*{0.5cm} & \text{in} \; \Omega\times(0,\tmax),\\
			v_t= \Delta v - \gamma u v, & \text{in} \; \Omega\times(0,\tmax),\\
			w_t= \Delta w-\delta u w+ \mu_3 w(1-w), & \text{in} \; \Omega\times(0,\tmax),\\
			u(x,0)=u_0(x), \quad v(x,0)= v_0(x), \quad w(x,0)= w_0(x), & x\in\overline{\Omega},
		\end{dcases}
	\end{equation}
	where the second and third equation are decoupled. The first two equations recall the consumption model studied in \cite{LankeitWang17} if we set $\gamma=k=1$.

	{\em Lemma \ref{lem:masslimBoundedness}} still holds by the parabolic maximum principle, although equation \eqref{eq:boundV} becomes
	\begin{equation*} 
		\lVert v(\cdot,t)\rVert_{\Lspace{\infty}{\Omega}} \leq M_3 := \lVert v_0\rVert_{\Lspace{\infty}{\Omega}}, \qquad \forall \ t\in (0,\tmax).
	\end{equation*}
	
	Regarding the boundedness of the gradients in Lemma \ref{lem:gradient}, after setting  $\alpha=\beta=0$, the only hypothesis required is then $\delta > \mu_3$; meaning that the rate of consumption of tumor cells by the lymphocytes has to be greater than the growth coefficient of the former. Such a requirement is consistent with biological intuition.
	
	The proof of Theorem \ref{theo:globalexist} relies on the ODI obtained for the functional $\Phi(t)$ in \eqref{eq:Phidef}, i.e., the proof of {\em Lemma \ref{lem:boundPhi}}. Without the problem parameters $\alpha$ and $\beta$, the systems of equations in \eqref{eq:systemK} and \eqref{eq:systemK1} have two parameters less to be chosen, namely $\Cr{e:dFgradvI}$ and $\Cr{e:dFgradvII}$. The proof of the theorem remains valid and so the conclusion.
	
	For the biological model, it could be not appealing decoupling the second and third equation in \eqref{eq:model} as \eqref{eq:modelWithoutAB} do. The lack of connection between $w$ and $v$ in \eqref{eq:modelWithoutAB} makes no sense, the latter being the chemical signal released by the former. However, as an extension of the model in \cite{LankeitWang17} for $k>1$, it is worthwhile to realize that large values of 
	\begin{math}
		\chi \lVert v_0 \rVert_{L^\infty(\Omega)}
	\end{math}
	are allowed without the need to impose further restrictions on the parameters, as in \eqref{eq:muLankeit} and \eqref{eq:muRaf}.
	
	The model in \cite{GnanasekaranCDFN24} does not admit global classical solutions for large values of	
	\begin{math}
		\lVert w_0 \rVert_{L^\infty(\Omega)}.
	\end{math}
	We can notice, as an advantage of the model proposed in this work, that the presence of a source term for the first equation in \eqref{eq:model} makes the existence of such solutions independent on 
	\begin{math}
		\lVert w_0 \rVert_{L^\infty(\Omega)}.
	\end{math}
	It would also be interesting to study whether the global solution stabilizes for large times as in \cite{LankeitWang17} and \cite{HuTao20}. After rewriting the first equation in \eqref{eq:model} as 
	\begin{math}
		u_t = \Delta u-\chi \nabla\cdot(u \nabla v) + \mu_1 u^k \left( 1 - \tfrac{\mu_2}{\mu_1} u \right),
	\end{math}
	we expect from the logistic term that $\tfrac{\mu_1}{\mu_2}$ signifies the carrying capacity of the lymphocytes (i.e., $u$). Moreover, aligned with the observations in \cite{HuTao20}, being the amount of lymphocytes bounded in time, they eventually eliminate both cancer cells (i.e., $w$) and the chemical signal (i.e., $v$). Therefore, we conjecture that 
	\begin{equation*}
		\left(  u(\cdot,t), v(\cdot,t), w(\cdot,t) \right) \rightarrow \left(\tfrac{\mu_1}{\mu_2}, 0, 0 \right) \quad \textbf{in}\;\; \Lspace{\infty}\Omega, \qquad 
		\textbf{as}\;\; t \rightarrow \infty.
	\end{equation*}
	
	%
	%
	
	\section*{Acknowledgments}
	
	The author is member of the Gruppo Nazionale per l’Analisi Matematica, la Probabilità e le loro Applicazioni (GNAMPA) of the Istituto Nazionale di Alta Matematica (INdAM) and acknowledges financial support by PNRR e.INS Ecosystem of Innovation for Next Generation Sardinia (CUP F53C22000430001, codice MUR ECS0000038).


\begin{thebibliography}{99}
		
		\bibitem{BaghaeiKelghati17}
		Baghaei, K. and Khelghati, A. (2017).
		\newblock Boundedness of classical solutions for a chemotaxis model with
		consumption of chemoattractant.
		\newblock {\em C. R. Math. Acad. Sci. Paris}, 355(6):633--639.
		
		\bibitem{BellomoBTW15}
		Bellomo, N., Bellouquid, A., Tao, Y., and Winkler, M. (2015).
		\newblock Toward a mathematical theory of {K}eller--{S}egel models of pattern
		formation in biological tissues.
		\newblock {\em Math. Models Methods Appl. Sci.}, 25(9):1663--1763.
		
		\bibitem{Evans10PDEs}
		Evans, L. (2010).
		\newblock {\em Partial Differential Equations}, volume~19 of {\em Graduate
			Series in Mathematics}.
		\newblock American Mathematical Society, 2nd edition.
		
		\bibitem{GnanasekaranCDFN24}
		Gnanasekaran, S., Columbu, A., {D\'iaz Fuentes}, R., and Nithyadevi, N. (2025).
		\newblock Global existence and lower bounds in a class of tumor-immune cell
		interactions chemotaxis systems.
		\newblock {\em Discrete Contin. Dyn. Syst. - S}, 18(6).
		
		\bibitem{HorstmannSurvey03}
		Horstmann, D. (2003).
		\newblock From 1970 until present: the {K}eller--{S}egel model in chemotaxis
		and its consequences. i.
		\newblock {\em Jahresber. Dtsch. Math.-Ver.}, 105(3):103--165.
		
		\bibitem{HorstmannSurvey04}
		Horstmann, D. (2004).
		\newblock From 1970 until present: the {K}eller--{S}egel model in chemotaxis
		and its consequences. ii.
		\newblock {\em Jahresber. Dtsch. Math.-Ver.}, 106(2):51--69.
		
		\bibitem{Horstmann2005}
		Horstmann, D. and Winkler, M. (2005).
		\newblock Boundedness vs. blow-up in a chemotaxis system.
		\newblock {\em J. Differential Equations}, 215(1):52--107.
		
		\bibitem{HuTao20}
		Hu, B. and Tao, Y. (2020).
		\newblock Critical mass of lymphocytes for the coexistence in a chemotaxis
		system modeling tumor–immune cell interactions.
		\newblock {\em Z. Angew. Math. Phys.}, 71(167).
		
		\bibitem{KellerSegel70prod}
		Keller, E. and Segel, L. (1970).
		\newblock Initiation of slime mold aggregation viewed as an instability.
		\newblock {\em J. Theor. Biol.}, 26(3):399--415.
		
		\bibitem{KellerSegel71cons}
		Keller, E. and Segel, L. (1971).
		\newblock Traveling bands of chemotactic bacteria: a theoretical analysis.
		\newblock {\em J. Theor. Biol.}, 30(2):235--248.
		
		\bibitem{LankeitWang17}
		Lankeit, J. and Wang, Y. (2017).
		\newblock Global existence, boundedness and stabilization in a high-dimensional
		chemotaxis system with consumption.
		\newblock {\em Discrete Contin. Dyn. Syst.}, 37(12):1078--0947.
		
		\bibitem{LankeitWinklerSurvey20}
		Lankeit, J. and Winkler, M. (2020).
		\newblock Facing low regularity in chemotaxis systems.
		\newblock {\em Jahresber. Dtsch. Math. Ver.}, 122:35--64.
		
		\bibitem{LankeitWinklerSurvey23}
		Lankeit, J. and Winkler, M. (2023).
		\newblock Depleting the signal: Analysis of chemotaxis-consumption models -- a
		survey.
		\newblock {\em Stud. Appl. Math.}, 151:1197:1229.
		
		\bibitem{LiLankeit16}
		Li, Y. and Lankeit, J. (2016).
		\newblock Boundedness in a chemotaxis-haptotaxis model with nonlinear
		diffusion.
		\newblock {\em Nonlinearity}, 29(5):1564--1595.
		
		\bibitem{Mahlbacher19}
		Mahlbacher, G., Reihmer, K., and Frieboes, H. (2019).
		\newblock Mathematical modeling of tumor-immune cell interactions.
		\newblock {\em J. Theor. Biol.}, 469:47--60.
		
		\bibitem{Painter19}
		Painter, K. (2019).
		\newblock Mathematical models for chemotaxis and their applications in
		self-organisation phenomena.
		\newblock {\em J. Theor. Biol.}, 481:162--182.
		
		\bibitem{Tao11}
		Tao, Y. (2011).
		\newblock Boundedness in a chemotaxis model with oxygen consumption by
		bacteria.
		\newblock {\em J. Math. Anal. Appl.}, 381(2):521--529.
		
		\bibitem{TaoWinkler12A}
		Tao, Y. and Winkler, M. (2012a).
		\newblock Boundedness in a quasilinear parabolic-parabolic {K}eller--{S}egel
		system with subcritical sensitivity.
		\newblock {\em J. Differential Equations}, 252(1):692--715.
		
		\bibitem{TaoWinkler12C}
		Tao, Y. and Winkler, M. (2012b).
		\newblock Eventual smoothness and stabilization of large-data solutions in a
		three-dimensional chemotaxis system with consumption of chemoattractant.
		\newblock {\em J. Differential Equations}, 252(3):2520--2543.
		
		\bibitem{TelloWrzosek16}
		Tello, J. and Wrzosek, D. (2016).
		\newblock Predator-prey model with diffusion and indirect prey-taxis.
		\newblock {\em Math. Models Methods Appl. Sci.}, 26(11):2129--2162.
		
		\bibitem{TelloWinkler07}
		Tello, J.~I. and Winkler, M. (2007).
		\newblock A chemotaxis system with logistic source.
		\newblock {\em Commun. Partial Differ. Equ.}, 32(6):849--877.
		
		\bibitem{WangWang20}
		Wang, J. and Wang, M. (2020).
		\newblock The dynamics of a predator-prey model with diffusion and indirect
		prey-taxis.
		\newblock {\em J. Dyn. Diff. Equat.}, 32:1291--1310.
		
		\bibitem{Wilkie17}
		Wilkie, K.~P. and Hahnfeldt, P. (2017).
		\newblock Modeling the dichotomy of the immune response to cancer: Cytotoxic
		effects and tumor-promoting inflammation.
		\newblock {\em Bull. Math. Biol.}, 79(6):1426--1448.
		
		\bibitem{Winkler10}
		Winkler, M. (2010).
		\newblock Aggregation vs. global diffusive behavior in the higher-dimensional
		{K}eller--{S}egel model.
		\newblock {\em J. Differential Equations}, 248(12):2889--2905.
		
		\bibitem{Winkler18}
		Winkler, M. (2018).
		\newblock Finite-time blow-up in low-dimensional {K}eller--{S}egel systems with
		logistic-type superlinear degradation.
		\newblock {\em Z. Angew. Math. Phys.}, 69(40).
		
	\end{thebibliography}

	\appendix
	
	\section{Computation of constants} \label{sec:contants}
	
	This section is devoted to showing the constants computed during the text. It supports the statement that we have computable constants in our results, except the constants related to the Gagliardo-Nirenberg inequality which remain unknown. 
	
	As notation we define
	\begin{equation*} 
		\cY{\varepsilon}{r} := \tfrac{r-1}{r} (\varepsilon r)^{\frac{1}{1-r}},
	\end{equation*}
	meaning the constant dependent on $\varepsilon$ in the Young inequality: $a b \leq \varepsilon a^r + \cY{\varepsilon}{r} b^{\frac{r}{r-1}}$, for real positive $a,b,\varepsilon$ and $r>1$. This way makes easier the use of the inequality in the choice of parameters $r$ and $\varepsilon$ in each case.
	
	In the proof of Lemma \ref{lem:ineq}, some parameters are left as degrees of freedom for the a priori inequalities later on. Along with them, some constants dependent on them are computed, for instance, in \eqref{ineq:uPgradv2}, \eqref{ineq:uPgradv2KO}. They are:
	\begin{equation} \label{eqA:ceuPgradv2}
		\Cr{ce:uPgradv2}(\Cr{t:uYoung},\Cr{t:gradvYoung}):= \left[\cY{\Cr{t:uYoung}}{\tfrac{p+k}{p}}\right]^{\frac{(p+1)k}{p(k-1)}}  \cY{\Cr{t:gradvYoung}}{\tfrac{(p+1)k}{p+k}} \measDom,
	\end{equation}
	\begin{equation} \label{eqA:ceuPgradv2KO}
		\Cr{ce:uPgradv2KO}(\Cr{t:uYoung}) := \cY{\Cr{t:uYoung}}{p+1},
	\end{equation}
	\begin{equation} \label{eqA:ceu2gradvDPMD}
		\Cr{ce:u2gradvDPMD}(\Cr{t:uYoung},\Cr{t:gradvYoung}):= \left[\cY{\Cr{t:uYoung}}{\tfrac{p+k}{2}}\right]^{\tfrac{(p+1)(p+k-2)}{2(k-1)}}
		\cY{\Cr{t:gradvYoung}}{\tfrac{(p+1)(p+k-2)}{(p-1)(p+k)}} \measDom,
	\end{equation}
	\begin{equation} \label{eqA:ceu2gradvDPMDKO}
		\Cr{ce:u2gradvDPMDKO}(\Cr{t:uYoung}) := \cY{\Cr{t:uYoung}}{\tfrac{p+1}{p-1}},
	\end{equation}
	\begin{equation} \label{eqA:cegradVgradW}
		\Cr{ce:gradVgradW}(\Cr{t:uYoung},\Cr{t:gradvYoung}):= \left[\cY{\Cr{t:uYoung}}{\tfrac{2(p+1)}{p-1}}\right]^{\tfrac32} \cY{\Cr{t:gradvYoung}}{3} \measDom,
	\end{equation}
	and
	\begin{equation} \label{eqA:cegradW}
		\Cr{ce:gradW}(\Cr{t:uYoung}) := 2\cGN C^{p(1-\theta)} \cY{\frac{\Cr{t:uYoung}}{2\cGN C^{p(1-\theta)}}}{\frac1\theta} + 2\cGN C^p.
	\end{equation}
	
	In Section \ref{sec:localExist}, it is defined
	\begin{equation} \label{eqA:cmassY}
		\Cr{c:massY} := \tfrac{\mu_1}{k+1} \left(\tfrac{2 \mu_1 k}{\mu_2 (k+1)} \right)^p \measDom.
	\end{equation}
	
	In Section \ref{sec:aprioriest}, during the proof of Lemma \ref{lem:gradient}, there are computed
	\begin{equation} \label{cons:indterm}
		\Cr{c:indterm} := \Cr{c:normsum}\measDom\left(1 + 2M_1 + \frac{\Cr{c:linODI}M_1}{\mu_2} \right)+  \Cr{c:normsum}\cY{\frac{1}{2}}{\frac{k+1}{2}}\measDom^{\frac{k+1}{k-1}} + \frac{\mu_1 \Cr{c:normsum}}{\mu_2}\cY{\frac{\mu_2}{2\mu_1}}{\frac{k+1}{k}}\measDom^{k+1}
	\end{equation}
	and
	\begin{equation} \label{cons:indtermU}
		\Cr{c:indtermU} := \Cr{c:normsum} \measDom + \left(\frac{\mu_1 \Cr{c:normsum}}{\mu_2} + \tfrac{\Cr{c:normsum}\Cr{c:linODI}}{\mu_2} + 2\Cr{c:normsum} \right) M_1.
	\end{equation}
	
	Next, for the proof of Lemma \ref{lem:dFu}, it is defined
	\begin{equation} \label{eqA:dFu}
		\Cr{c:dFu} := \mu_1 p\  \cY{\frac{\mu_2}{2\mu_1}}{\frac{p+k}{p+k-1}} \measDom.
	\end{equation}
	Following, Lemma \ref{lem:dFgradv} relies on the constant
	\begin{equation} \label{eqA:dFgradv}
		\Cr{c:tempII}(\Cr{e:dFgradvI},\Cr{e:dFgradvII},\eta):= \Cr{c:LankeitineqU}(\eta) + (2\alpha p)^\frac32 \Cr{ce:gradVgradW}\left(\tfrac{\Cr{e:dFgradvI}}{2\alpha p},\Cr{e:dFgradvII}\right).
	\end{equation}
	Finally, Lemma \ref{lem:dFgradw} obtain the constant
	\begin{equation} \label{eqA:dFgradw}
		\Cr{c:tempIII}(\Cr{e:ineqGradW},\eta) := \Cr{c:LankeitineqU}(\eta) +  4p\mu_3 \cGN C^p \left(\cY{\frac{\Cr{e:ineqGradW}}{4p\mu_3 \cGN C^{p(1-\theta)}}}{\frac1\theta} C^{-p\theta} + 1\right),
	\end{equation}
	for the constants $\theta= \frac{p-1}{p-1+\frac2n}$, $C$ from \eqref{eq:boundednessgradv}, and $\cGN$ from the Gagliardo-Nirenberg inequality as cited in \eqref{ineq:gradW}.
	
\end{document}